\catcode`\@=11
\def\magnification{\afterassignment\m@g\count@}
\def\m@g{\mag\count@}
 \documentclass[10pt,a4paper]{article}
\usepackage{amsmath,amsfonts}
\usepackage{eurosym}
\usepackage{graphicx,color}
\graphicspath{{/Users/christian/Documents/Dogbe/POSTSCRIPT/}}
\DeclareGraphicsExtensions{,.ps,.eps}

\newcommand{\sect}[1]{\setcounter{equation}{0}\section{#1}}
 
\newcommand{\seqn}{\begin{equation}}
\newcommand{\eeqn}{\end{equation}}
\newcommand{\seqna}{\begin{eqnarray}}
\newcommand{\eeqna}{\end{eqnarray}}

\newtheorem{theo}{Theorem}[section]
 \newtheorem{lem}[theo]{Lemma}
\newtheorem{pro}[theo]{Proposition}
\newtheorem{cor}[theo]{Corollary}
\newtheorem{hypo}[theo]{Hypothesis}

\let\bb\mathbb
\def \R {{\relax{\bb R}}}

\def \D {\relax{\bb D}}

\def\S {\relax{\bb S}}

\def \d {\partial }

\def \a {\alpha}

\def \eps {\epsilon }

\def\a {\alpha}
\def\o {\omega}

\def\la{\langle }
\def \ra {\rangle}\def \s {\sigma}
\def \frac#1#2{{\textstyle{#1 \over #2}}}\def\o{\omega }
\def \eps {\epsilon }
 
 \def \d {\partial }

\begin{document}

\title{Spherical Harmonics Expansion of the Vlasov-Poisson  initial boundary value problem\\ }

\author{Christian   {DOGBE}\thanks {e-mail: 
dogbe@math.unicaen.fr}\\
Laboratoire de Math\'ematiques Nicolas Oresme, \\ CNRS, UMR 6139\\
Universit\'e de Caen\\  
BP 5186, 14032 Caen Cedex
} 

\maketitle
\date 

\begin{abstract} 
We derive and analyze the `SHE' (Spherical Harmonics Expansion) type 
system of equations  coupled in energy.  We also show that
diffusive behavior occurs on long time and distance scales and we
determine the diffusion tensor.  The analysis is based on the
governing kinetic equations  arising in electron transport in
semiconductors.   
\end{abstract}

	

Ê

\smallskip
 \section{Introduction}
 
This paper presents an extension of a previous work by P. Degond
and S. Mancini  \cite{DM} in which a diffusion model describing
the evolution of an electron gas confined between two parallel
planes by a strong magnetic field is derived. Electrons colliding against the plane are supposed to be reemitted following a combination of diffusive and specular elastic laws. When the distance between these planes goes to zero, the distribution function (solution of a Vlasov-Poisson equations) is proved to converge to the solution of the macroscopic model, the so-called SHE model  which is a diffusion model for the energy distribution function.   This situation arises in a certain kind of ion propellers for satellites. 
SHE is an acronym for Spherical Harmonics Expansion coming from its earlier derivation by physicists \cite{DSR}. These models appeared to be a good compromise between a very accurate description of physical phenomena by kinetic models and less numerically expensive macroscopic models such as Drift-Diffusion or Energy-Transport models. Whereas kinetic equations deal with functions of seven variables (six variables in the phase and time) and classical macroscopic models deal  with functions of four variables (three in the position, space  and time), SHE models introduce an intermediate one dimensional variable replacing the velocity and which appears as the kinetic energy associated to the velocity of a particle. This makes of SHE models quite reliable models at a rather low  numerical cost. The SHE  model in  literature has been derived from the Boltzmann equation first by P. Dmitruk, A. Saul, and L. Reyna \cite{DSR}.

\quad
The mathematical theory of the diffusion approximation started
with the seminal papers \cite{BLP}, \cite{BSS} in the context of
neutron transport after the formal theory had been set up by
Hilbert, Chapman, Enskog and co-workers (see for instance
\cite{CH}). Their approaches were later extended by many authors.
In this work, we are interested in a physical situation where the observation length scale is large compared to the mean free path while the observation time is large compared to the characteristic time evolution of the particle.

\smallskip
\quad
The mathematical study of the vanishing mean free path limit and
of diffusion approximation is by now  a classical problem with
applications in various fields of physics. We refer among others
to \cite{BBP}, \cite{LK} and for recent   \cite{D2} and for SHE
models  \cite{BD}, \cite{DM} and the references therein for details as the physical background concerning these  models.
Concerning   SHE  models this study   has been derived from the Boltzmann equation first by  N. Ben Abdallah and P. Degond \cite{BD} under the assumption that the dominant scattering mechanisms are elastic collisions.

\quad
The present work has been inspired by \cite{BBP},   \cite{D2},  \cite{DM}, \cite{DLMM}, but differs from them in the nature of collisions of electrons with
the wall and the presence of a magnetic field.
In \cite{DM} and \cite{DLMM},  elastically diffusive collisions at the plates is considered: the particles are reemitted with the same energy as their incident one, and with a random velocity direction. 
The force field is the gradient of a smooth potential function which is assumed given, independent of time and of the mean free path, and amounts to assume vary  over the macroscopic scale only.
As a consequence, the large time behavior of the distribution function is given by an energy distribution function depending on the parallel components relative to the boundary of the domain   and satisfying a diffusion equation in both position and energy. Here, we treat the self-consistent problem, the electric potential satisfying the Poisson equation.
Hence, in contrast of  \cite{DM} and \cite{DLMM}, the presence of a quadratic term in Vlasov equation gives rise to some singular term which   adds additional technical difficulties. The existence and regularity of the solution to the asymptotic model which is now constituted of the diffusion equation coupled to  the Poisson equation, are not immediate, nevertheless verified. Our main contribution in this paper is to give a rigorous proof of this convergence.

\quad
In \cite{BBP} the authors showed that collisions with a boundary
can drive  a system of neutral particles towards a diffusion
regime and  the diffusivity is infinite. 
It has been shown in \cite{D2} that a logarithmic time re-scaling restores
a finite diffusivity.  Nevertheless  the divergence which appears
in \cite{BBP} and waived  in \cite{D2} does not appear in
our work because of the presence of a strong magnetic field
directed transversally to the plates. Particle motion then does
not occur along straight lines, but rather along the plates whose
axes are parallel to the magnetic field lines.  

We note also, that some results in the direction of this paper have been obtained, without the boundary conditions, and  where collisions are taken into account through the non linear Pauli operator by T. Goudon et {\it al} in \cite{GM}.

\sect{Setting up the problem and the main result}\label{SETPROBLEM}

\quad
Let us recall   the physical background carefully. Our starting point is a
Vlasov  equation for the electron distribution. The  microscopic
model  describes the evolution of particles limited between two
parallel plates. In the region of this space, the particles are
submitted to a force field. The physical space
where the electrons evolve is $\Omega\times \R^3$ where $\Omega
=[0,1]\times \R^2$.  The position vector is denoted by   $X=(x,y,z)$ and   we split   $X$ into its perpendicular  $ x\in[0,1]$ and parallel components
$\underline \xi =(y,z)\in\R^2$ relative to the boundary $\Gamma$.
We denote by $v=(v_x, v_y, v_z)\in \R^3$ the position and the
velocity vectors of an electron between the planes. We  decompose   
$v= (v_x, \underline v)$, where $v_x$ is the velocity component
parallel to the $x$-axis and $\underline v=(v_y, v_z)\in \R^2$ is the
component parallel to the plates.  The electrons are subject to a
magnetic field transverse to the plates $(B(\underline \xi), 0, 0)$
depending upon
$\underline
\xi$  and  to a potential force field  parallel  to the plates,
$(0, E_y(\underline \xi), E_z(\underline \xi))$   depending only
upon $\underline \xi$ and  satisfying   $\underline E
=-\nabla_{\underline\xi}\phi$,  where $\nabla_{\underline\xi}$
denotes the 2-dimensional gradient with respect to $\underline
\xi$ and $\phi$ is the potential.  (cf. Fig.1). 
Therefore, they are supposed to move between the plates
according to a collisional transport equation.
\seqna
\nonumber
 \includegraphics[scale=1.]{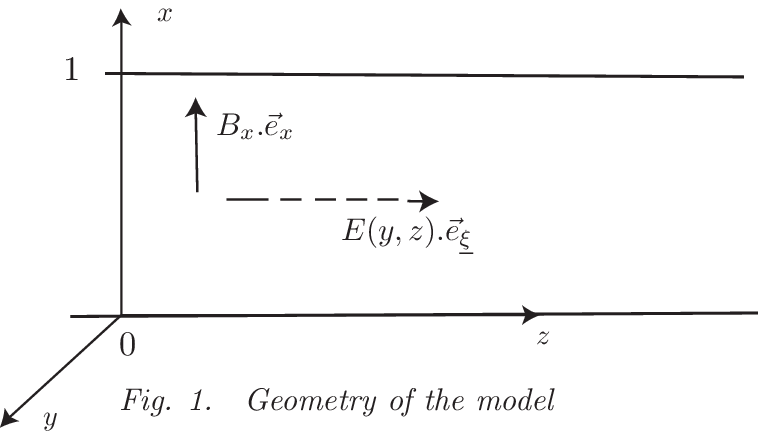} 
\eeqna
Here
$f(X,v,t)dXdv$ is understood to give the number of electrons that
occupy any infinitesimal volume  $dXdv$ at the point $(X,v)$ which
is given as the solution of the boundary value problem:
\seqna
\label{Boltz1}
{\d_t f} +
\left(\underline v \cdot\nabla_{\underline \xi} + {\underline
E} \cdot \nabla_{\underline v}\right) f  + v_x{\d f \over \d x} -
(\underline v\times B)\cdot \nabla_{\underline v} f  =0
\eeqna
where  $(\underline v\times B) = (0, v_zB, -v_yB)$, while the
electrostatic potential $\phi$ solves the Poisson equation
\seqna
\label{POIS1}
\underline E =-\nabla_{\underline \xi}\phi,\qquad -\Delta_{\underline
\xi} \phi  =\int_{0}^1\!\!\!\int_{\R^3} f dxdv - C(\underline \xi). 
\eeqna
We shall assume that $C(\underline\xi)$ is a given
(nonnegative) function which is regular, say $C^\infty(\Omega)$.

\quad
Let us specify the proper boundary conditions to be
considered. We introduce the set $\Theta =\Omega\times \R^3$, and
consider its boundary $\Gamma = \gamma\times \R^3$, where  $\gamma =\
\{0, 1\}\times \R^2$; we denote $\R_\pm ^3 =\{v\in \R^3\,\, :\,\,
\pm v_x>0\}$ 
and the following {\it incoming} and {\it
outgoing} subsets of $\Gamma$ (respectively representing incoming
and outgoing particles to the domain $\Omega$):
 \seqna
\label{GAMMABORD} 
\nonumber
\Gamma_- =\left(\{0\}\times \R^3\times \R_+^3\right)\cup
\left(\{1\}\times \R^3\times \R_-^3\right)
\eeqna
where for instance, $\left(\{0\}\times \R^3\times \R_+^3\right)$
represents electrons entering the region $0<x<1$, through the boundary
plane $x=0$. The boundary  $\Gamma_+$ is obtained by reversing the
inequalities. We introduce the traces (or boundary values) of
$f$ on
$\Gamma\times \R^3$ according to:
 \seqna \label{Bord0} \gamma(f)
=\left.f\right|_{\{x=0,1\}},\qquad \gamma^{\pm}(f) =
\left.f\right|_{\{x=0,1,\,\,\,\pm v_x>0\}},
 \eeqna
$\gamma^{+}(f)$ is the distributional function of the particles exiting
the domain $\Omega$ at the boundary $\Gamma$, while $\gamma^{-}(f)$ is
that of incoming particles. We suppose, finally that it is a function
of outgoing trace through an operator ${\cal K}$ which expresses  the
interaction of the particles:
\seqna
\label{Bord1}
\gamma^-(f) ={\cal K}(\gamma^+(f)).
\eeqna
Introducing the operator 
${\cal B}(\gamma (f)) =  \gamma^-(f)- {\cal
K}(\gamma^+(f))$, we have ${\cal B}(\gamma(f))=0$.
The boundary condition for $\phi$ is:
\seqna
\label{BOUNDPHI}
\lim_{|\underline \xi|\to\infty}\phi(\underline\xi, t)=0,\quad
\mbox{a.e.}\,\,\,\, t>0.
\eeqna
The homogeneous boundary condition for $\phi$ means that the system
of electrons is in equilibrium at infinity.  We refer to \cite{CIP}
for further properties of these boundary conditions and physical
interpretations. In this work some important relations valid for
these boundary conditions are derived.

\smallskip
Our model is based on the assumption  that the distance  
between the  plates  is small compared to the characteristic length
of the horizontal motion. Since we are looking for a diffusion
process, this suggests to rescale the longitudinal coordinate   
according to  
$x'= \alpha x$, where $\alpha =O(1)$ is  the small parameter. But
then, the number of collisions between a typical particle and the
plates per unit of unscaled time is large: if we assume  that the
collisions between the particles and the plates are a purely isotropic
process, there is no  reason to expect that, at  a large scale, the
particles would follow a horizontal drift. Hence, in order to
observe a horizontal motion  at a large scale, it is logical to
rescale the time variable as $t'=\alpha^2 t$. After rescaling,
(dropping  the primes for the sake of clarity), the equation
(\ref{Boltz1})-(\ref{POIS1}) is recast as: 
 \seqna
\label{MainModel1} 
&& \alpha {\d_t f^\alpha} \,\,+\,\,
\left(\underline v \cdot\nabla_{\underline \xi}  +
{\underline E}^\alpha \cdot\nabla_{\underline v}\right)
f^\alpha + {1\over \alpha} \left(v_x {\d  \over \d x}  -
(\underline v \times B)\cdot
\nabla_{\underline v}\right)  f^\alpha =0,\\
\label{MainModel2}
&& {\underline E}^\alpha =  -\nabla_{\underline \xi} \phi^\alpha,\quad
-\Delta_{\underline \xi} \phi^\alpha =    \int_{0}^1\!\!\!\int_{\R^3} f^\alpha dxdv -C(\underline \xi).
\eeqna
It is clear that these scalings do not induce any modification in the
boundary condition, namely, $f^\alpha$ still satisfies:
\seqna
\label{Bord3} 
\gamma^-(f^\alpha) ={\cal K}(\gamma^+(f^\alpha)).
\eeqna 
Finally,  we prescribe an initial condition which is
compatible with the expected asymptotic dynamics. This initial
condition is 
\seqna 
\label{Boltz-init} f (X, v, t=0) &=&
f_I(X,v)\quad \forall (X,v)\in \Omega\times \R^3.
\eeqna
Now,  since we want to base our approximation of diffusion
in position-energy space, it is  relevant to interpret the velocity
of particle in terms of its energy. For that purpose, let us
introduce the spherical coordinates in velocity space:
\seqna
\label{COAREA1}
\left\{\begin{array}{cc}
\omega =v/|v|\in
\S^2 \, (\hbox{angular variable});\,\,\,\,\,\\
\eps = |v|^2/2,  (\hbox{energy variable});\,\,\,\qquad \quad\\
v= |v|\omega =\sqrt{2\eps}\omega= N(\eps)\omega. \qquad\qquad
\end{array} \right.
\eeqna
We identify $v$ with the pair $(\varepsilon, \omega)$. Hence, for
$C^1$ function $\eps : \R^3 \to \R$, and any integrable
function   $\varphi \in
\R^3$  we have (coarea formula \cite{F}):
\seqna
\label{COAREA2}
\int_{v\in \R^3}\varphi(v)dv &=& \int_0^\infty
\int_{\S^2}\varphi (\sqrt{2\eps} \, \omega)(2\eps )^{1/2}{d\eps\,
d\omega}\\
\nonumber
&=& \int_{\eps>0}\int_{\omega\in\S^2} \varphi(\eps,
\omega) N(\eps)d\eps  d\omega
\eeqna
where $d\omega$ denotes the normalized Euclidean measure on
$\S^2$, so $dv = r^2drd\omega$, with  $r=|v| =\sqrt{2\eps }$, and
$dr ={d\varepsilon\over \sqrt{2 \eps}}$.
In the coarea formula (\ref{COAREA2}), $N(\eps)$ represents the
energy-density of states.   

\smallskip
\smallskip
The main contribution of this paper is to investigate the limit
$\alpha\to 0$ of (\ref{MainModel1})-(\ref{Boltz-init}). We show that
the limit
$f^0$ as $\a\to0$ is a function $F(\underline \xi, \varepsilon, t)$
of the longitudinal coordinate $\underline \xi$, of the energy
$\varepsilon =|v|^2/2$ and of the time which obeys a diffusion
equation in the position-energy space.  In order to simplify the exposition, from now on, we will denote by $(P)$ the problem (\ref{MainModel1}), (\ref{MainModel2}),  (\ref{Bord3}), (\ref{Boltz-init}).

\smallskip
\quad 
We  are now ready to state our main result.
\begin{theo}
\label{MAINTHEOREM}
\begin{enumerate}
\item
As $\a$ goes to zero, the solution $f^\alpha$ to  the problem 
$(P)$ formally converges to an
equilibrium state  $F(\underline \xi,|v|^2/2, t)$, solution to the
following  SHE model: 
 \seqna 
\label{CONT}
4\pi\sqrt{2\eps} {\d F\over \d t}   &+&
\left(\nabla_{\underline\xi}- {\underline E}{\d\over
\d\eps}\right)\cdot {\underline J} =0,\\
\label{CONT-CURR}
{\underline J} ({\underline\xi}, \varepsilon, t)&=& -
\D({\underline\xi},\varepsilon) \left(\nabla_{\underline\xi}- 
{\underline E}{\d\over \d\eps} \right)F, \\
 \label{MainModel2Theo}
{\underline E}  &=& -\nabla \phi,\quad -\Delta_{\underline \xi}=
\int_0^1\!\!\!\int_{\R^3} f(\underline \xi,v,t)  dxdv-C(\underline \xi), \\
\label{NEWINIT}
F({\underline\xi}, \varepsilon, t=0)
&=&F_I({{\underline\xi}}, \varepsilon), \\
\label{CURR}
 \underline J ({\underline\xi}, \varepsilon =0, t) &=&0,
\eeqna
where $F_I$ is a suitable initial condition, in the domain
$(\xi,\eps)\in \R^2\times (0,\infty)$. The  `diffusivity tensor'
$\D$ is   given by:
\seqna
\label{COEFFDIFFUSION1}
\D({\underline\xi},\varepsilon)
=(2\eps)^{3/2}\int_0^1\!\!\!\int_{\S^2}\underline
\chi(x,\omega; {\underline\xi}, \eps)\otimes \underline \omega
dxd\omega,
\eeqna
where ${\underline \omega} =(\omega_y, \omega_z)$,
$\,\,\,\underline \chi= (\chi_y, \chi_z),\,\,\, $
$\,\,{\underline \chi}\otimes
\underline\omega$ is the tensor product
$(\chi_i\omega_j)_{i,j\in\{y,z\}}$ and  $\chi_i(x,\omega;
{\underline\xi},\frac{|v|^2}2)$  is a solution of the problem
\begin{equation}
\label{EQUATIONHOMOLOGIQUE}
\left\{\begin{array}{cc}
\displaystyle -v_x{\d \chi_i\over \d x} + (\underline v\times B)\cdot
\nabla_{\underline v} \chi_i= \omega_i&\quad \mbox{in}\,\,\,
\Theta  \\
\\ \displaystyle \gamma^+(\chi_i) ={\cal K}^*(\gamma^-(\chi_i))
\qquad\qquad \quad&\quad\mbox{on}\,\,
\,\Gamma,
\end{array}\right.
\end{equation}
 satisfying $\displaystyle \int_0^1\!\!\!\int_{\S^2}\underline \chi
dxd\omega=0$, for all $(\underline \xi, \varepsilon)$, where  ${\cal K}^*$ is the operator adjoint of ${\cal K}$.

\item
Under hypothesis \ref{HYPOKERNEL}, \ref{HYPOCOUCHEINIT} (to be
specified later on) and
\ref{Hypotrans}, there exists a solution $(f^\a , \underline E^\a)$
in $L^\infty(0, T; L^2(\Theta))\cap L^2(0, T; L_{loc}^2(\R_{\underline
\xi}^2))$ for any $T\in \R^+$. $f^\a$ converges to $f^0$ in the weak
star topology of $L^\infty(0, T; L^2(\Theta))$ for any $T>0$, where
$f^0(t, \underline \xi, \varepsilon)= F(t, \underline \xi,
\varepsilon)$ and $F(t, \underline \xi, \varepsilon)$ is the weak
solution of the problem  (\ref{MainModel1})-(\ref{Boltz-init}).

\end{enumerate}
\end{theo}

\smallskip
\smallskip
The result is very close to \cite{DM}. Nevertheless, the proof will be
different, though we shall use a lot of results developed in
\cite{DM}.
Existence results of Cauchy problem for Vlasov-Poisson system are
now well known (see \cite{Pa}) in dimension three without external
potential and the  boundary   value problem  for Vlasov-Poisson was
studied in \cite{Gy}. The theory of global weak solutions is due to
\cite{Ar}. The detailed mathematical study of
$(P)$ in the case of a constant
magnetic field and a given  electrostatic potential, which varies
only over the macroscopic scale will be found in \cite{DM}. We shall
not dwell on the existence and uniqueness of a solution for the
Cauchy problem $(P)$  which can be
done by means of Leray-Schauder's fixed point Theorem; we shall
solely give a material which  makes it  possible to obtain this solution and
focus on  the  establishment of the limit model.   The investigation
of this limit when $\a\to 0$ proceeds in two steps. The first one
consists in showing that $f^\a$ formally converges to a function of
$f(x, \eps, t)$ only, solution of diffusion  `SHE' problem. Then we
derive the continuity and current equation (\ref{CONT}),
(\ref{CONT-CURR}).
 The second one   corresponds to  a rigorous convergence proof; we 
show that weak solutions of Vlasov-Poisson equations converge
weakly (in an $L^2$ sense) towards weak solutions of the `new' SHE
model.   To achieve these goals, two methods can be developed: the
Hilbert expansion method \cite{BSS} and the moment method  \cite{D2}.
We shall choose the latter because it involves more straightforward
computations.

The outline of the paper is as follows.  Section
\ref{BOUNDOPERATION} is devoted to the preliminary materials
regarding the functional setting of our problem and important
properties of the boundary collision operator. We establish some
mathematical results on this operator.  In particular, we prove that
the flux of particles is conserved at the boundary and there is an
entropy dissipation (that is, a Darroz\`es-Guiraud-like inequality in gas surface interaction \cite{DG}, \cite{CIP}).
This will be enough to allows us to derive the formal asymptotic
limit which is achieved in Section  \ref{FORMAL}. This asymptotic limit  appears as a singular perturbation problem for
Vlasov-Poisson equations. We solve this problem and derive
our SHE model as the corresponding limit equation.  A necessary step
for the rigorous derivation is to establish the existence of a solution
for the kinetic problem; this will be done in Section
\ref{EXISTENCESOLUTION}. Finally  Section
\ref{RIGOROUS} is dedicated to the rigorous proof of convergence
itself.

\sect{The boundary operator: assumptions and
properties.}\label{BOUNDOPERATION}

In the sequel, we consider an expression of ${\cal K}$ as follows
\seqna
\label{Bord2}
 {\cal K} \phi  (v) =\int_{\{\omega'\in
\S^2,\,\, (\underline\xi, v)\in \Gamma_+\}} K(\underline \xi ,
|v|^2/2;
\omega'\to
\omega)
\phi(|v|\omega')|\omega_x'|d\omega',\qquad \forall v\in \R^3,
\eeqna
where $\omega =(\omega_x, \omega_y, \omega_z)={v\over|v|},\,\,\,\,\,
|\omega|=1$  is the decomposition of $v$ into spherical coordinates,
$\S^2=\{ \omega \in \R^3,\,\, |\omega|=1\}$ is the unit sphere.
The operator ${\cal K}$ maps the outgoing trace
$\gamma^{+}(f)$ to the incoming one $\gamma^{-}(f)$ .
The expression (\ref{Bord2}) models an elastic bounce against the
planes with a random deflection of the velocity.
Note that $K$  is an integral kernel which describes the
reflection law of the velocity direction. The quantity $K(\underline
\xi, |v|^2/2;
\omega'\to \omega)|\omega_x|d\omega$ is the probability for
a particle hitting the planes at point $\underline \xi$ with velocity
$v' =|v|\omega'$ to be reflected with the same $|v|$ and velocity
direction $\omega$ in the solid angle $d\omega$. Bearing in   mind
that ${\cal K}={\cal K}(\underline \xi, |v|^2/2)$ operates on the
angular variable $\omega$ only while ${\underline \xi}$ and $|v|$
are mere parameters.

 \smallskip 
We give  the adapted functional setting for the study of the collision
operator. We set   ${\cal S}_\pm(x),\,\, x=0,1$ to be the following
half-spheres:
\seqna
{\cal S}_+(0) ={\cal S}_-(1) =\{\omega\in \S^2,\,\,\,
\omega_x<0\},\quad {\cal S}_-(0) ={\cal S}_+(1) =\{\omega\in
\S^2,\,\,\, \omega_x>0\}.
\eeqna
We introduce the domain ${\cal S} =[0,1]\times
 \S^2$ with its associated {\it incoming} and {\it outgoing}
boundaries defined by:
\seqna
{\cal S}_-= \left({\cal S}_-(0) \times \{0\}\right)\cup
\left({\cal S}_-(1) \times \{1\}\right),\qquad
{\cal S}_+= \left({\cal S}_+(0) \times \{0\}\right)\cup
\left({\cal S}_+(1) \times \{1\}\right).
\eeqna
We denote by $L^2({\cal S}_\pm)$ the space of square
integrable functions on ${\cal S}_\pm$ with respect to the
measure
$|\omega_x|d\omega$; by $(f,g)_\Theta,\,\,\,\,
(f,g)_{\Gamma^\pm}$  the inner products on
$L^2(\Theta)$ and on $L^2(\Gamma^\pm)$ respectively defined by:
\seqna
\nonumber
(f,g)_\Theta =\int_\Theta fgd\theta,\qquad
(f,g)_{\Gamma^\pm} =\int_{\Gamma^\pm}fg|v_x|d\Gamma
\eeqna
where $d\theta= dxd\underline\xi dv$ is the volume element in phase
space, and $ d\Gamma =\sum_{x=0,1} d\underline{\xi}dv$ is
the surface element and by $(f,g)_{\cal S},\,\,
(f,g)_{{\cal S}_{\pm}}$ the inner products on $L^2({\cal S})$ and
$L^2({\cal S}_\pm)$   defined analogously:
\seqna
\nonumber
(f,g)_{\cal
S}&=&\int_0^1\!\!\!\int_{\S^2}(fg)(x,\omega)dxd\omega,\\
\nonumber
(f,g)_{{\cal
S}_{\pm}}&=&\int_0^1\!\!\!\int_{{\cal
S}_{\pm}(0)}(fg)(x,\omega)|\omega_x|dxd\omega + \int_0^1\!\!\!\int_{{\cal
S}_{\pm}(1)}(fg)(x,\omega)|\omega_x|dxd\omega.
\eeqna
and $|f|_{{\cal S}},\,\, |f|_{{\cal S}_\pm}$ the associated
norms. We shall denote also $\Theta' = \R^2\times \R^{+*}$ with
$d\theta'=d\underline x d\varepsilon$ its volume element.
$\forall f = f(\o), \,\,\, \la f\ra$ will denote the angular average
of $f$ on the sphere $\S^2$ and with respect to $x$ i.e.
\seqna
\nonumber
\la f \ra  ={1\over 4\pi} \int_0^1\!\!\!\int_{\S^2} f(x, \o) dxd\o.
\eeqna
We denote by ${\cal C}^\pm =\{ f\in L^2(\R^3_\pm):\,\, \,\, f(\eps,
\o)\,\,
\mbox{is constant with respect to}\,\,\, \o \}$.
We define the operator $Q^\pm$ as the orthogonal projection 
(for inner product $(\cdot, \cdot)_{{\cal S}_+}$) of $L^2({\cal
S}_+)$ on the space  ${\cal C}^\pm$, i.e. 
\seqna
Q^\pm  f(x,\omega) =
{1\over \pi}\int_{{\cal S}_{\pm}(x)}  f(\omega) |\omega_x|d\omega,
\qquad
\o\in {\cal S}_\pm,\quad x\in \{0, 1\}, 
\eeqna
and the operator  $P^\pm$,  as the orthogonal complement of $Q^\pm :
\,\,\, P^\pm = I-Q^\pm$.

\bigskip
 \quad 
We shall  list the required  properties of the reflection
operator $K$. They are summarized in the following
\begin{hypo}
\label{HYPOKERNEL}
(i) Flux conservation:
\seqna
\label{con}
\int_{{\cal S}_-(x)} K(\omega'\to \omega)|\omega_x|d\omega =1,
 \eeqna
$\mbox{for almost all}\,\,\,\,(\omega,\omega')\in {\cal
S}_-(x)\times{\cal S}_+(x),\,\,\,\, x=0,1$.

\smallskip 
\noindent
\item(ii) Reciprocity principle:
\seqna
\label{rec}
K(\omega'\to \omega) =K(-\omega\to -\omega'),\qquad \forall
(\omega, \omega')
\in {\cal S}_-(x)\times{\cal
S}_+(x),\,\,\,\, x=0,1.
\eeqna
$ (iii)\,\,  \mbox{Positivity}:
\label{iii}
K(\omega'\to \omega)>0,\quad
\mbox{for almost all}\,\,\,\,(\omega,\omega')\in {\cal
S}_-(x)\times{\cal S}_+(x),\,\,\,\, x=0,1$.
 
\noindent
(iv) The operator ${\cal K}(\underline \xi, \varepsilon)$ is a compact operator
from $L^2({\cal S}_+)$  onto $L^2({\cal S}_-)$.
\end{hypo}

\smallskip															
The relation (\ref{con}) expresses  the conservation of the normal flux
of particles at the boundary.

\noindent
The reciprocity relation (\ref{rec}) is a macroscopic effect
of the time reversibility of elementary particle-surface
interactions.
As a direct consequence of the flux conservation and reciprocity
relation, by the change of $\omega$ into $-\omega$ in (\ref{con})
and the use of (\ref{rec}), we get the following `normalization'
identity:
 \seqna \label{norm} 
\int_{{\cal S}_+(x)} K (\omega'\to
\omega)|\omega'_x|d\omega'=1,\qquad x=0,1. 
\eeqna 
The assumption
(iv) can be viewed as  a  regularity hypothesis for the integral
kernel $K$. The
adjoint of ${\cal K}(\underline \xi, \varepsilon)$ is obviously
given by 
\seqna
 \label{ADJOINT-B} 
 {\cal K}^*(\phi) (x, \o)  
=\int_{{\cal S}_-(x)} K(x, \o\to \o')\,\phi(x, \o')\,|\o_x'|\,
d\omega',\qquad \o\in {\cal S}_+(x),\quad x=0,1,
\eeqna
 and is also a compact operator.

\subsection{Flux conservation}

Our first objective is to understand what  the effects of the
boundary condition (\ref{Bord2})  are on the fluxes of particles.
For this, we need to  define the incoming and outgoing normal fluxes
associated to distribution function $f$.
\seqna
\nonumber
J_x^- (\underline \xi, \varepsilon, t) &=&  \int_{\o\in {\cal S}_-(x)}
f(\underline \xi, x=0, |v|\o|) |v||\o_x|d\o\\
\nonumber
J_x^+ (\underline \xi, \varepsilon, t) &=&
\int_{\o'\in {\cal S}_+(x)}f(\underline \xi, x=1, |v|\o'|)
|v||\o_x'|\, d\o'. 
\eeqna
We begin with the following the required identity:
\begin{lem}
\label{LEMMAFLUX}
Assume ${\cal B}(\gamma(f^\a))=0$. Then the normal flux
is conserved:
\seqna
\label{FLUX0}
J_x^-=J_x^+.
\eeqna
\end{lem}

\noindent
{\bf Proof.}
The incoming normal flux $J_x^-$ is obtained  by integrating the
left-hand side of the equation
\seqna
\gamma^-(f^\a) =  {\cal K}(\gamma^+(f^\a)).
\eeqna
Indeed, the coarea formula implies
\seqna
\nonumber
J_x^- (\varepsilon) &=&  \int_{\o\in {\cal S}_-(x)}
\gamma^-(f^\a)(\varepsilon \o) \sqrt{2\varepsilon}\, |\o_x| d\o
\\
\nonumber
J_x^+ (\varepsilon) &=&
\int_{\o\in {\cal S}_+(x)}
\gamma^-(f^\a)(\varepsilon \o) \sqrt{2\varepsilon}\, |\o_x| d\o.
\eeqna
Hypothesis (\ref{con}) and Fubini's Theorem lead to
\seqna
\label{FLUX1}
J_x^- (\varepsilon) &=&  \int_{\o\in {\cal S}_-(x)}
\left(\int_{\o'\in {\cal S}_+(x)} K(\o'\to \o')f(|v|\o'|)|\o_x'|\,
d\o'\right)  |v||\o_x|\,d\o\\
\nonumber
&=& \int_{\o'\in {\cal S}_+(x)}f(|v|\o'|) |v||\o_x'|\,
d\o' = J_x^+(\varepsilon).
\eeqna
Finally, it is easy to check that the flux particles through
$\Gamma_\pm$  vanished.

\smallskip
\noindent 
The hypothesis \ref{HYPOKERNEL} is crucial for the validity
of Theorem \ref{MAINTHEOREM}. From this, the  following inequality
bears similarities with Darroz\`es-Guiraud inequality in gas
dynamics, which we can deduce from Jensen's inequality:
\begin{lem}
\label{LEMMADARROGURAU}
Let $\gamma^+(f)\in L^2({\cal S}_+)$ and $\gamma^-(f)= {\cal
K}(\gamma^+(f))$. Then
\seqna
\label{DARROGURAU}
\int_{{\cal S}_-(x)}  |\gamma^-(f (x,\omega))|^2|\omega_x|\,d\omega
\leq \int_{{\cal S}_+(x)}   |\gamma^+(f(x,\omega))|^2\omega'_x|
d\omega'.
\eeqna
\end{lem}
 
\smallskip
\noindent
{\bf Proof.}  Using Cauchy-Schwarz inequality and normalization
identity (\ref{norm}), we have for $x=0,1$:
\seqna 
\int_{{\cal S}_-(x)}  |\gamma^-(f (x,\omega))|^2|\omega_x|\,d\omega
&=& \int_{{\cal S}_-(x)}    |{\cal K}(\gamma^+(f(x,\omega)))|^2
|\omega_x| d\omega \\
\nonumber
&\leq& \int_{{\cal S}_-(x)}  \left|\int_{{\cal S}_+(x)}  K(x,
\o'\to \o)\,\gamma^+(f)(x, \o')\,|\o_x'|
\,d\o'\right|^2|\o_x|\,d\o\\
\nonumber
&\leq&  \int_{{\cal S}_+(x)}   |\gamma^+(f(x,\omega))|^2\omega'_x|
d\omega'.
\eeqna
Therefore, the operator norm of ${\cal K}$ in the space
${\cal L}( L^2({\cal S}_+), L^2({\cal S}_+))$ of bounded operator  
from $L^2({\cal S}_-)$ to $L^2({\cal S}_+)$ is less or equal one:  $\|{\cal
K}\|_{{\cal L}(L^2({\cal S}_+))} \leq1$. Furthermore, from
(\ref{norm}), if $\varphi$ is constant over ${\cal S_+}$, then ${\cal
K}\varphi$ is constant over ${\cal S_-}$. Hence, $\|{\cal K}\|_{{\cal
L}(L^2({\cal S}_+))} =1$.

\smallskip
\quad
We introduce the specular reflection operator ${\cal J}$
which operates from $L^2(S_+)$ to  $L^2(S_-)$ according to ${\cal J}
\varphi (\omega) =\varphi(\omega^*)$, with $\omega_* =(-\omega_x,
\omega_y, 
\omega_z)$. Its adjoint  ${\cal J}^*$  is defined from
$\phi\in L^2({\cal S}_-)$ to $ L^2({\cal S}_+)$  and is also the mirror
reflection operator, and we have the equalities
$\,\,{\cal J}^*{\cal J}  =I_{{\cal S}_+} = {\cal J}{\cal J}^*
=I_{{\cal S}_-}$. Therefore the operator ${\cal K}{\cal J}^*$ and its
adjoint ${\cal K}^*{\cal J}$  operate on $L^2({\cal S}_-)$ while
${\cal K}^*{\cal J}$ and ${\cal K}{\cal J}^*$ operate on $L^2({\cal
S}_+)$.

\smallskip
In the derivation of the SHE model, we are interested in the
characterization of the equilibria of operator $I- {\cal K} {\cal
J}^*$ and  $I- {\cal K}^* {\cal J}$ in $L^2({\cal S}_\pm)$. 
By this hypothesis (iv), the operators $I- {\cal K} {\cal J}^*$, $I-
{\cal K}^* {\cal J}$ are  Fredholm operators. Using Krein-Rutman's
Theorem and Fredholm's theory  \cite{KR}, we have the following
lemma  which can be easily adapted from that of \cite{DM}, so the
proof is omitted:

\begin{lem}
\label{Nullspace}
(i) The null-spaces $N(I- {\cal J}{\cal K}^*)$ and
$N(I- {\cal J}^*{\cal K})$ are spanned by the constant
functions on ${\cal S}_-$ and ${\cal S}_+$ respectively.

\noindent
\item(ii) ${\cal K}(x, \varepsilon)$ is of norm 1, i.e. $\|{\cal K}(x, \varepsilon)\varphi\|_{{\cal L}(L^2({\cal S}_+))}
\leq\|\varphi\|_{{\cal L}(L^2({\cal S}_+))}$, for all
$\varphi\in L^2({\cal S}_+)$.

\noindent
\item(iii) The range $R(I- {\cal K}^*{\cal J})$ is such that
$R(I- {\cal K}^*{\cal J})= N(I- {\cal J}^*{\cal K})^\bot$.
Equivalently, the equation $(I- {\cal K}^*{\cal J}) f=g$ has a
solution $f$ if and only if  $\displaystyle\int_{{\cal
S}_+}g(\omega)|\omega_x|d\omega =0$.  Then, the solution
$f$ is unique under the additional constraint
$\displaystyle\int_{{\cal S}_+}f(\omega)|\omega_x|d\omega =0$.
\end{lem}
Property (ii) is   reminiscent of Darroz\`es-Guiraud inequality for
boundary conditions in gas dynamics. This property expresses that ${\cal K}(\underline \xi, \varepsilon)$ has a `good' diffusion behavior.
From (\ref{Bord2}) , ${\cal K}\varphi(\omega)$ for
$\omega\in {\cal S}_-$ appears as a convex mean value of
$\varphi(\omega')$ over ${\cal S}_+$. 

\begin{lem}
\label{ORTHOPROJECT}
The operator ${\cal K}$ satisfies
\seqna
\label{ORTHOPROJECT1}
{\cal K} Q^+ = Q^-{\cal K} ={\cal J}Q^+ ={\cal J} Q^-,\qquad {\cal
K}P^+ = P^-{\cal K}.
\eeqna
\end{lem}

\smallskip
The proof of this lemma is a straightforward adaptation from
Ref. \cite{DM}. 
 
\smallskip
We close this section by giving   the following assumption:
\begin{hypo}
\label{HYPOPROJECT}
There exists  $k_{0}<1$ such that, for $|v|\in \R^+,\,\,
\underline \xi\in \R^2$

\smallskip
\centerline{$\|{\cal K} P^+\|_{{\cal L}(L^2({\cal S}_+),
L^2({\cal S}_-))}\leq k_{0}<1$.}
\end{hypo}
This assumption comes from the elementary operator theory. Indeed, for every $\underline \xi\in \R^2$ and $|v|\in \R^+$, there exists $k(\underline \xi, |v|)$ such that ${\cal K} P^+\|\leq k(\underline \xi,|v|<1$.  Obviously, the hypothesis \ref{HYPOPROJECT} is satisfied in the case of the isotropic scattering.

\sect{Formal derivation of the macroscopic model}\label{FORMAL}

As a  further purpose, we  consider a sequence $(f^\a)_\a$ of solutions to the problem  $(P)$. Our goal is to study the asymptotic
behavior of solutions $(P)$  as $\alpha\to0$. In this section, we assume that $f^\alpha \to f$ as $\alpha \to 0$ in a smooth way (in the sense where we will precise later).
This approach will formally give rise to a set  of
diffusion equations to be viewed as a SHE system. In fact, this formal
derivation will not be completed in this section since an `auxiliary
problem' - the resolution of which will be postponed to Section 
\ref{RIGOROUS}- will arise during this study.

\smallskip
We divide the formal asymptotic into several steps. First, we prove
that the limit $f$ is a function of total energy. Then we prove that
the limit satisfies  a continuity equation. An auxiliary  problem will
be considered which allows us to finally derive the current equation.

\quad
The following change of variables will be useful in the
remainder of the paper.
 Since the velocity $v$ satisfies the relation.
\seqna
\label{PARAM1}
|v| = \sqrt{v_x^2+v_y^2+ v_z^2}, \qquad\,\,\, 
  \omega_j ={v_j\over |v|},\qquad j=x,y, z,\,\,\,
\eeqna
we  can  parameterize the sphere $\S^2$ in the
direction $\underline\omega= (\omega_y, \omega_z)$, and $\omega_x$ is
defined by
$\omega_x=\sigma\sqrt{1-\omega_y^2-\omega_z^2},\,\,\,\,
\sigma =\pm1$. Thus the sphere $\S^2$ is given by two local maps
${\{(\s, \omega_y, \omega_z),\,\,\, \sigma =\pm1\}}$.

\smallskip
Hence, we have 
\begin{lem}
\label{PARAM2}
Let the change of spherical coordinate be:
\seqna
(v_x, v_y, v_z) \longmapsto (v,\omega_x, \omega_y, \omega_z)
\eeqna
where
$\displaystyle \omega_x = \mbox{sgn}(v_x) \sqrt{1-\omega_y^2
-\omega_z^2},\quad \omega_ y = {v_y\over
|v|}, \quad   \omega_ z = {v_z\over  |v|}.$
Then, for $\mbox{sgn}(v_x)=\pm1$, we have
\seqna
\label{PARAM3}
\begin{array}{cc}
\begin{array}{cc}
\displaystyle {\d f\over  \d v_x} ={\d f\over \d
|v|}\cdot \omega_x - {\d f\over \d \omega_y}\cdot
{\omega_x\omega_y\over |v|} - {\d f\over \d \omega_z}
\cdot{\omega_x\omega_z\over |v|}\\
\displaystyle {\d f\over  \d v_y} ={\d f\over \d
|v|}\cdot \omega_y - {\d f\over \d \omega_y}\cdot
{1-  \omega_y^2\over |v|} - {\d f\over \d \omega_z}
\cdot{\omega_y\omega_z\over |v|}\\
\displaystyle {\d f\over  \d v_z} ={\d f\over \d
|v|}\cdot \omega_z - {\d f\over \d \omega_z}\cdot
{\omega_z\omega_z\over |v|} - {\d f\over \d \omega_z}
\cdot{1- \omega_z^2\over |v|}\\
\end{array}
\end{array}
\eeqna
\end{lem}
\smallskip\noindent
{\bf Proof.} - For $\mbox{sgn}(v_x) =\pm1$, we have
\seqna
\nonumber
{\d f\over \d v_j} ={ {\d f\over \d v}{\d v\over \d v_j} } +{{\d f
\over \d \omega_y} {\d \omega_y \over \d v_j}} +{{\d f
\over \d \omega_y} {\d \omega_z \over \d v_j}}.
\eeqna
Therefore, for $j=x, y, z$, we get
\seqna
\label{PARAM4}
{\d |v| \over \d v_j} ={v_j \over \sqrt{v_x^2 +v_y^2
+v_z^2}}=\omega_j.
\eeqna
Using  (\ref{PARAM4}) and the relation  $v_j =|v| \omega_j$, we
find
\seqna
\nonumber
{\d \omega_y \over \d v_y} ={1-\omega_y^2\over |v|},
\qquad\quad  {\d
\omega_z \over \d v_y} =-{\omega_y\omega_z\over
|v|};\\
\nonumber
{\d \omega_y \over \d v_z} =-{\omega_y\omega_z \over |v|},
\qquad\quad  {\d \omega_z \over \d v_z} ={1- \omega_z^2\over
|v|};\\
\nonumber
{\d \omega_y \over \d v_x} =-{\omega_y\omega_x \over |v|},
\qquad\quad {\d \omega_z \over \d v_x} =-{\omega_z\omega_x\over
|v|},
\eeqna
which proves the result.

\subsection{The limit is a function of the energy}

Taken    formally the limit $\a\to 0$
in (\ref{MainModel1})-(\ref{Bord3}) shows that $f$ is the solution to
equations:
\seqna
\label{FunctEner4E1}
\left(v_x {\d  \over \d x}  - (\underline v \times B)\cdot
\nabla_{\underline v} \right) f & =& 0, \\
\label{FunctEner4E2}
\gamma^-(f) &=&{\cal K}(\gamma^+(f)).
\eeqna  
We show the
\begin{lem}
\label{ENERGYFUNCTION} 
  The solution of (\ref{FunctEner4E1})-(\ref{FunctEner4E2}) are
functions of total energy only: 
\seqna
\label{FunctEner2}
f(x,v, t) = F(\underline \xi,  \frac12|v|^2 , t).
\eeqna
\end{lem}
 
\smallskip\noindent
{\bf Proof.}  
 The fact that  the limit of $f^\alpha$ is a function of the
energy  only is an easy consequence of  Lemma \ref{Nullspace}. 
According to Lemma  \ref{PARAM1},     
Eq.  (\ref{FunctEner4E1}) is equal to \seqna
\label{FunctEner5} {\d f \over \d x}(x, \omega)   +
{B(\underline\xi)\over |v|\omega_x} {\d f\over \d
\omega}(e_x\times \omega) =0 \eeqna where $\displaystyle {\d
f\over \d \omega}$ denotes the derivatives of $f$ with respect to
$\omega\in \S^2$ of degree 1, $a \times b$ stands for the cross
product of two vectors $a$ and  $b$ and $e_x \times \omega$ denotes
the tangent  vector to $\S^2$ in $x$-axis. We recall that
$\underline \xi$ and $v$ are mere parameters in problem
(\ref{FunctEner4E1}) and will be omitted in the remainder of the
proof. For $f(x, \omega) $ we introduce the following change of
variables: 
\seqna 
\label{FunctEner6}
 \left\{\begin{matrix}
\displaystyle \theta ={B({\underline \xi})\over |v|\omega_x} =
{B({\underline \xi})\over \sqrt{2\eps}\omega}, \quad \omega_x
=\s\sqrt{1-\omega_y^2-\omega_z^2},\quad \s\in\{-1,
1\}\\
\omega_y^* =\omega_y\cos(\theta  x)+ \omega_z\sin(\theta
x),\qquad \qquad \qquad\qquad\qquad\qquad\quad\,\,\,\,\\
\omega_z^* =   \omega_y\sin(\theta x)+ \omega_z\cos(\theta
x).\qquad \qquad \qquad\qquad\qquad\qquad\quad\,\,\,\,
\end{matrix} \right.
\eeqna
Since $\underline \omega =(\omega_y, \omega_z)$, with
(\ref{FunctEner6}) we can write, $\underline \omega^*=
(\omega_y^*, \omega_z^*) ={\cal R}^+_{x, \s}(\underline \omega)$,
where
\seqna
\label{FunctEner7}
{\cal R}_{(x,\sigma)}^+(\underline\omega) =
\left(\begin{matrix}
\cos \theta x&-\sin \theta x\\
\sin \theta x&\cos \theta  x
\end{matrix}\right)
\eeqna 
is a rotation of $\underline \omega$ around the $x-$axis by
an  angle $\theta x$. Hence the function associated to this
transformation is written as \seqna f^*(x,\sigma,\underline
\omega) = f(x,\sigma,{\cal R}_{(x,\sigma)}^+(\underline\omega)).
\eeqna
 Therefore, if $f$ is solution of (\ref{MainModel1}),  $f$
satisfies the problem: \seqna \label{FunctEner8} |v|\omega_x  {\d
f^* \over \d x}  =0, \qquad \gamma^-(f^*) ={\cal K}(\gamma^+(f^*)). 
\eeqna
 Integrating the first equation of
(\ref{FunctEner8}) with respect to $x\in[0,1]$, yields
 \seqna
f^*(x, \sigma, \underline\omega) =f^*(0, \sigma, \underline\omega)
=f^*(1,\sigma,\underline\omega), \eeqna or, in terms of the
rotation, ${\cal R}_{(x,\sigma)}^-(\underline\omega)$: \seqna
f(x,\sigma,\underline\omega) =f^*(0,\sigma, {\cal
R}_{(x,\sigma)}^-(\underline\omega)) =f(1,\sigma, {\cal
R}_{(1-x,\sigma)}^+(\underline\omega)). 
\eeqna 
This is  equivalently written
 \seqna 
\label{FunctEner11}
\gamma^+(f) ={\cal K}^*\gamma^-(f).
\eeqna
Inserting (\ref{FunctEner11}) into the second equation of
(\ref{FunctEner4E2}) leads to
\seqna
(I-{\cal K}{\cal J}^*)(\gamma^-(f))=0
\eeqna
which, by virtue of  Lemma \ref{Nullspace}, implies that
$\gamma^-(f) =F(\frac12|v|^2)$.

\subsection{The continuity equation}

We are now concerned with the derivation of the continuity equation
(\ref{CONT}).  Before, we introduce the following macroscopic
quantities. For finite $\a>0$, we define the average density
$F^\alpha$ and current ${\underline J}^\a(\underline\xi, \eps,t) =
(J_y^\a, J_z^\a)$ by
\seqna
\label{DENSITY}
F^\a(\underline \xi, \varepsilon, t) &=& {1\over 4\pi}
\int_0^1\!\!\!\int_{\S^2} f^\a (x, \underline \xi, |v|, \o, t) dx d\o,\\
\nonumber
J^\a(\underline \xi, \varepsilon, t) &=&  {|v|\over
\a}\int_0^1\!\!\!\int_{\S^2} \underline v f^\a(x, \underline \xi,|v|,
t)dxd\o : ={1\over \a}N(\varepsilon)\la v f^\a (\varepsilon,
\cdot)\ra\\
\label{CURRENT}
&=& {2\eps\over \a}\int_0^1\int_{\S^2}\underline \o f^\a (x,
\underline
\xi, \varepsilon, \o, t) dx d\o.
\eeqna
According to reference \cite{DM}, we have,
\begin{lem}
\label{CURRENTDIFFERENTIAL1}
Let $\varphi(x, v)$ be a $C^1$ function. 
Then we have
\seqna
\label{CURRENTDIFFERENTIAL2}
{\d J_\varphi\over \d \varepsilon} = \sqrt{2\varepsilon}
\int_0^1\!\!\!\int_{\S^2} (\nabla_{\underline v}\varphi)(x,
\sqrt{2\varepsilon} \omega)dxd\omega,\qquad
J_\varphi(\varepsilon) = 2\varepsilon
\int_0^1\!\!\!\int_{\S^2}\underline \omega \varphi(x,
\sqrt{2\varepsilon} \omega)dxd\omega.
\eeqna
\end{lem}

\smallskip
We prove:
\begin{lem}
\label{CONTINUITY1} 
We temporarily admit that $J^\a\to J$ as $\a\to0$.
Then $F$ and $J$ satisfy the continuity equation (\ref{CONT}).
\end{lem}

\smallskip\noindent
{\bf Proof.}  
Since from Lemma \ref{ENERGYFUNCTION}, $f$ is  a function of energy
only (in velocity space), it is quite natural to expect that $f^\alpha$
be mainly given by  a function of $\eps$ only  at order  O  in $\alpha$,
for instance the macroscopic quantity $\la f^\alpha\ra$. We  expect
that there exists a function  $g^\alpha (\underline \xi, v, t)$ such
that:
\seqna
\nonumber
f^\alpha (X, v, t) = F^\alpha (\underline \xi, \eps, t) + \alpha
g^\alpha (\underline \xi, v, t).
\eeqna
We integrate (\ref{MainModel1}) with respect to $x$ and
on a sphere of constant energy. We get
\seqna
\label{FORMAL4}
\nonumber
\int_0^1\!\!\!\!\int_{\S^2}\left(\underline v \cdot\nabla_{\underline
\xi}   - {\underline E}^\alpha \cdot\nabla_{\underline v}
\right)(F^\alpha + \alpha g^\alpha)dxd\o &=&
\int_0^1\!\!\!\!\int_{\S^2} {\underline v}d\o\cdot {\widetilde
\nabla}F^\alpha + \alpha \nabla_{\underline
\xi}\cdot \int_0^1\!\!\!\!\int_{\S^2}{\underline v}g^\alpha dxd\o\\
\nonumber
&& -\alpha \underline E^\alpha \cdot
\int_0^1\!\!\!\!\int_{\S^2}\nabla_{\underline v} g^\alpha dxd\o\\
\nonumber
&=&{4\pi\alpha \over N(\eps)}\nabla_{\underline \xi}\cdot J^\alpha 
- \alpha \underline E^\alpha \cdot
\int_0^1\!\!\!\!\int_{\S^2}\nabla_{\underline v} g^\alpha dxd\o
\eeqna
with ${\widetilde \nabla} =\nabla_{\underline \xi} +\underline E \d/
\d \eps$ and for $\displaystyle \int_{\S^2} vd\o=0$.
Moreover, for any test function $\psi(\eps)$, a straightforward
computation yields:
\seqna
\label{FORMAL6}
\nonumber
 \int_{\eps>0}\psi \int_0^1\!\!\!\int_{\S^2}\nabla_{\underline v}
g^\alpha dxd\o\, N(\eps)d\eps = 4\pi \int_{\eps>0}\psi {\d J^\alpha
\over \d\eps} d\eps
\eeqna
which finally proves
\seqna
\label{FORMAL7}
{1\over 4\pi \alpha} \int  (\underline v \cdot\nabla_{\underline \xi}  -
{\underline E}^\alpha \cdot\nabla_{\underline v})f^\alpha dxd\o
N(\eps)
 = {\widetilde \nabla} \cdot J^\alpha.
\eeqna

\subsection{The current equation} 

\smallskip
Next, we have to find a current equation giving a relation
between $J^\alpha$ and the density $F^\alpha$. We prove:
 
\begin{lem}
\label{CURRENT2} The current equation (\ref{CONT-CURR}) is satisfied. 
\end{lem}

\smallskip
\noindent
{\bf Proof.}   Using Green's formula  (\ref{LEMMAGREEN-1}) we get
\seqna
\label{CURR0}
S^\alpha &:=& \int_0^1\!\!\!\int_{\R^3} \left( v_x{\d  \over \d x}
 - (v\times B)\cdot \nabla_{\underline v} f^\alpha \right) \chi dxdv\\
\nonumber
&=&\left(\int_{\R^3}v_xf^\alpha {\underline \chi}
dv\right)\left|_{x=1}\right. 
-\left(\int_{\R^3}v_xf^\alpha {\underline \chi}
dv\right)\left|_{x=0}\right.
+ \int_0^1\!\!\!\int_{\R^3} f^\alpha\left(v_x{\d f^\alpha\over \d x}
 -  (v\times B)\cdot \nabla_{\underline v}\right){\underline
\chi}dxdv.
\eeqna
Now, using (\ref{Bord3}) and the second equation of
(\ref{EQUATIONHOMOLOGIQUE}) we get 
\seqna
 \label{CURR1}
 \nonumber
\left.\left(\int_{\R^3}v_xf^\alpha {\underline  \chi}
dv\right)\right|_{x=1} &=& |v|^2 \left(\int_{{\cal S}^+}
\gamma^+(f^\alpha) \gamma^+ (\underline \chi)|\omega_x|d\omega   -
\int_{{\cal S}_-} \gamma^-(f^\alpha) \gamma^- (\underline \chi)
|\omega_x|d\omega\right)\\
\nonumber
 &=& |v|^2 \left (\int_{{\cal S}_+}  \gamma^+(f^\alpha)
\gamma^+ (\underline \chi)|\omega_x|d\omega  -\int_{{\cal S}_-}
{\cal K}( \gamma^+(f^\alpha)) \gamma^- (\underline
\chi)|\omega_x|d\omega\right)\\
&=& |v|^2  \int_{{\cal S}_+} \gamma^+(f^\alpha) 
\left (\gamma^+ (\underline \chi)|\omega_x|d\omega  -  {\cal
K}^*(\gamma^- (\underline \chi)\right)|\omega_x|d\omega=0. 
\eeqna
Similarly,  we can prove the relation
for $\left.\left(\int_{\R^3}v_xf^\alpha {\underline  \chi}
dv\right)\right|_{x=0}$. With the first equation of
(\ref{EQUATIONHOMOLOGIQUE}) we get
\seqna
\label{CURR2} 
 \int_0^1\!\!\!\int_{\R^3} f^\alpha\left(v_x{\d
\over \d x}
- (v\times B)\cdot \nabla_{\underline v}\right){\underline
\chi} dxdv &=& 2\eps \int_0^1\!\!\!\int_{\R^3} \underline \omega f^\alpha  
\underline \chi dxd\omega =  \alpha J^\alpha. 
\eeqna
Therefore, inserting (\ref{CURR1}) and  (\ref{CURR2}) into
(\ref{CURR0}), we deduce that $S^\alpha =\alpha J^\alpha$. 
 This justifies  the definition of ${\underline
\chi}$.
From (\ref{MainModel1}) and (\ref{CURR0}), we deduce:
\seqna
{\underline J}^\alpha = -\alpha \int_0^1\!\!\!\int_{\R^3} \d_t f^\alpha
{\underline \chi} dxdv - \int_0^1\!\!\!\int_{\R^3}
\left(\underline v \cdot\nabla_{\underline \xi} -
\nabla_{\underline \xi}\phi\cdot \nabla_{\underline v}\right) 
 f^\alpha {\underline \chi} dxdv.
\eeqna
Taking the limit $\alpha \to 0$ and using (\ref{FunctEner2}), we obtain
\seqna
\label{CURRENTSATISFIED}
{\underline J} =  -\int_0^1\!\!\!\int_{\R^3} 
\left(\underline v
\cdot\nabla_{\underline
\xi} - \nabla_{\underline \xi}\phi\cdot \nabla_{\underline v}\right) 
F {\underline \chi} dxdv,
\eeqna
where $\underline J =\lim {\underline J}^\alpha, \, F=\lim F^\alpha$.
Taking into account the relation
\seqna
\nonumber
\left(\underline v \cdot\nabla_{\underline
\xi} - \nabla_{\underline \xi}\phi\cdot \nabla_{\underline v}\right)
F^\alpha ={\underline v}  \left(\nabla_{\underline\xi} -
\nabla_{\underline\xi}\phi{\d\over
\d_\eps}\right)F^\alpha   
\eeqna   
we get
\seqna
\label{CCURRENTT}
{\underline J}= -\left(\int_0^1\!\!\!\int_{\S^2}\underline
\chi(x,\omega; \xi, \eps)\otimes \underline \omega dxd\omega\right)
\left(\nabla_{\underline\xi} - \nabla_{\underline\xi}\phi{\d\over
\d_\eps}\right) F, 
\eeqna
which leads to equations (\ref{CONT-CURR}).  
Finally, collecting  (\ref{CCURRENTT})  and
(\ref{FORMAL6}), inserting in (\ref{FunctEner4E1}), multiplying by
$\sqrt{2\varepsilon}$  and taking the limit $\alpha \to 0$ lead and second to 
(\ref{CONT}), (\ref{COEFFDIFFUSION1}).

It only remains to prove the existence of $\underline \chi$ which
will be achieved in section \ref{RIGOROUS}. 
The  relevance  of the model
(\ref{CONT})-(\ref{CONT-CURR}) relies on the positivity
of the coefficient $\D$. We postpone this task to  Section \ref{POSDIFFUSE}.

\sect{Existence of the solutions to the microscopic
problem}\label{EXISTENCESOLUTION}

\quad In this section, we establish the existence of solutions to the
microscopic problem (\ref{MainModel1}), (\ref{Bord3}). To avoid the
treatment of initial layers when we pass to the limit $\alpha\to0$,
we  impose a compatibility condition on the initial data $f_I$:
\begin{hypo}
\label{HYPOCOUCHEINIT}
There exists  a smooth function $F_I$ such that
$f_I(x,\underline \xi, v)=F_I(\underline\xi, |v|^2/2)$
 and that $f_I$ satisfies:  
$f_I\in L^2(\Theta),\,\, (\underline v
\cdot\nabla_{\underline \xi}  - \nabla_{\underline \xi}\phi
\cdot \nabla_{\underline v})f_I\in L^2(\Theta)$.
\end{hypo}


\bigskip
We define  the transport operator by:
 \seqna
\label{transp}
{\cal A}^\alpha f = \left(\underline v \cdot\nabla_{\underline
\xi} - \nabla_{\underline \xi}\phi\cdot \nabla_{\underline v}\right) 
 f + {1\over \alpha}\left( v_x{\d\over \d x}  - (\underline v
\times B)\cdot \nabla_{\underline v} \right)f 
\eeqna
on the  domain 
\seqna
\nonumber
D({\cal A}^\alpha) =\{ f\in L^2(\Theta),\quad {\cal A}^\alpha f \in
L^2(\Theta),\quad\gamma^+ (f)\in L^2(\Gamma_+),\quad  \gamma^-(f)
={\cal K} \gamma^+ (f)\}. 
\eeqna
We shall denote   by ${\cal A}$ the bare differential operator
(\ref{transp}) when no indication of the domain is needed. Following
\cite{B}  we define the space 
\seqna
H({\cal A}^\alpha) &=&\{ f\in L^2(\Theta),\quad {\cal A}^\alpha f \in
L^2(\Theta)\}.
\eeqna
According
to \cite{B}, it is well known that the regularity $f\in H({\cal
A}^\alpha)$ is not sufficient to guarantee the integrability of 
$|\gamma^+(f)|_{L^2(\Gamma_+)}^2$ and $|\gamma^-(f)|_{L^2(\Gamma_-)}^2$ over the
boundary. But if one of these traces is integrable, the other one is
also integrable. Therefore, we define the following space 
\seqna
\nonumber
H^\alpha_0 ({\cal A}) &=&\{ f\in  H^\alpha({\cal A}), \,\,
\gamma^-(f)\in L^2(\Gamma_-) \} =\{ f\in  H^\alpha({\cal A}),
\,\,\gamma^+(f)\in L^2(\Gamma_+)\},
\eeqna
endowed with the graph norm  and family of semi-norms, for $R>0$
\seqna
\label{GRAPHNORM}
|f|^2_{\atop H^\alpha({\cal A}^\alpha)} =|f|^2_{\atop L^2(\Theta)}
+|{\cal A}^\alpha f|^2_{\atop L^2(\Theta)},\qquad |f|^2_{\Gamma_\pm,R} =\int_{\Gamma_\pm,|v|\leq R} |v_x||f|^2
d\Gamma.
\eeqna
Note that hypothesis \ref{HYPOCOUCHEINIT} implies that $f_I\in
D({\cal A}^\alpha)$ for all $\alpha>0$. 
For the solvability of the  transport equation, we require the following
hypothesis:
\begin{hypo}
\label{Hypotrans}
\item{(i)} The first and second derivatives of the potential   belong  to the Sobolev space    $W^{1,\infty}$ (in other words,  $\nabla_{\underline \xi} \phi$ is bounded and globally lipschitz over $\R^2$). 
\item{(ii)}  The magnetic field $B=  B(\underline \xi)\in C^1\cap
W^{1,\infty}(\R_{\underline \xi}^2)$.
\end{hypo}

We give some results concerning  the trace  operators which are
necessary in our setting. Since we will need to handle nonlinear
functions of $f$ and their traces on the boundary, we have to study
for which nonlinear functions of $f$ traces can be defined.   
The problem of existence of a trace   is fundamental for the Cauchy problem  with boundary conditions. This problem was investigate by many authors, such   C. Bardos  \cite {B}.  M. Cessenat in \cite {Ce} studied it and applied it to neutron transport equation,  also S. Ukai in \cite{U} for the free transport operator,   E. Beals et {\it al} in \cite{BP}  for the abstract time-dependent linear kinetic equations and recently  S. Mischler  in \cite{M}     for Vlasov-Boltzmann equation.

\smallskip
\quad
According to Ref.  \cite{B}, \cite{BP}  and   hypothesis \ref{Hypotrans}
 the following Green's  type identity can be easily deduced.

 \begin{lem} 
\label{LEMMAGREEN}
For $f$ and $g$ in $H_0({\cal A}^\alpha)$  compactly supported with
respect to $v$, we have:
\seqna
\label{LEMMAGREEN-1}
({\cal A}^\alpha f, g)_{\Theta} + (f, {\cal A}^\alpha g)_{\Theta}
={1\over \alpha} \left((\gamma^+(f), \gamma^+(g))_{\Gamma_+} -
(\gamma^-(f), \gamma^-(g))_{\Gamma_-}\right).
\eeqna
\end{lem}

\smallskip
From that lemma we will deduce several estimates for the operator \ref{transp}  and their traces on the boundary.

\subsection{Resolution of an approximate kinetic problem}

\quad 
Applying Leray-Schauder's fixed point Theorem to the transport
operator ${\cal A}^\a$ with domain $D({\cal A}^\a)$ would be enough
to prove the existence for the kinetic problem. Nevertheless, since
the operator ${\cal K}$  acts like the identity on ${\cal
C}^\pm$, we are lacking some estimates on the traces. Following
\cite{DM}, we introduce an approximate problem. Let 
$\eta>0$ be a small approximating parameter; we approach the boundary
operator  by
\seqna
\label{PERTUBEDOPERATOR}
{\cal K}_\eta = P^+{\cal K} +{1\over 1+ \eta} {\cal J} Q^+.
\eeqna
We consider the associated initial value problem
\seqna
\label{EXISTENCESOLUT1}
&&\alpha \d_t f_\eta^\alpha +{\cal A}f_\eta^\alpha =0,\qquad
f_\eta^\alpha|_{t=0} =F_\eta,\\
\label{EXISTENCESOLUT2}
&&{\underline E}^\alpha = -\nabla_{\underline \xi}
\phi^\alpha, \quad -\Delta_{\underline \xi} \phi^\alpha = \int_{\R^3}
f_\eta^\a\,dv -C(\underline \xi),\\
\label{EXISTENCESOLUTBord} 
&&\gamma^-(f_\eta^\alpha) ={\cal K}(\gamma^+(f_\eta^\alpha)) 
\eeqna 
with domain $D({\cal A}_\eta^\alpha) =\{f\in H({\cal A}^\alpha),\quad 
\gamma^+(f)\in L^2(\Gamma_+),\quad \gamma^-(f) ={\cal K}_\eta
(\gamma^+(f))\}$.
We denote by   ${\cal A}_\eta^\alpha$ the transport operator ${\cal
A}$  on the  domain $D({\cal A}_\eta^\alpha)$.  
Thanks to the Lemma \ref{ORTHOPROJECT}, we have 
\seqna
\|{\cal K}_\eta P^+\|_{{\cal L}(L^2({\cal S}_-), L^2({\cal S}_+))}<1,
\qquad \forall\, \eta>0.
\eeqna

 \bigskip
We can now   easily adapt Proposition 4.1
and Lemmas 4.2 and 4.3 in Ref.  \cite{DM}. Namely we have:
\begin{pro}
\label{EXISTENCESOLO}
(i) For all $\eta>0$, and for all $F_\eta\in D({\cal A}_\eta^\alpha)$,
there exists a unique function $f_\eta^\alpha \in {\cal C}([0, T];
D({\cal A}_\eta^\alpha))\cap {\cal C}^1([0, T]; L^2(\Theta))$, that
solves (\ref{EXISTENCESOLUT1})-(\ref{EXISTENCESOLUT2}).

\item{(ii)}
We have also $|f_\eta^\a|_{L^2(\Theta)}\leq
|F_\eta|_{L^2(\Theta)},\quad |\a \d_t f_\eta^\a|_{L^2(\Theta)}
=\|{\cal A}^\a f_\eta^\a|_{L^2(\Theta)}\leq
\|{\cal A}^\a F_\eta^\a|_{L^2(\Theta)}$.

\item{(iii)} Let $F_I$ be as in hypothesis \ref{HYPOCOUCHEINIT}.
There exists a sequence $(F_\eta)_{\eta>0}$ such that $F_\eta\in
D({\cal A}_\eta^\alpha)$ and $ F_\eta \to F_I,\quad {\cal A} F_\eta
\to {\cal A} F_I\quad \mbox{in}\,\,\, L^2(\Theta)\,\, \mbox{weak
star}$, as $\eta\to0$.
\end{pro}

\subsection{Traces estimates and existence for the kinetic problem}

\quad
 In this section, we assume that boundary conditions of the form 
(\ref{Bord3}) are satisfied and prove the existence of the traces on
the boundary under suitable assumptions.
 We first establish that the projection
$P^+$ of the trace at the boundary of a function of $D({\cal A}^\a)$
is controlled by the graph norm. This is the key needed to prove the convergence of a solution to the problem  $(P)$.

\smallskip
We  prove

\begin{lem} 
\label{LEMMATRACE0}
We assume that $\phi$ satisfies hypothesis \ref{Hypotrans}.
If $f\in D({\cal A}_\eta^\alpha)$, then there exists a constant $C>0$
such that
\seqna
\label{PRO}
|P^-\gamma^-(f)|^2_{\atop L^2\left(\Gamma_-\right)}\leq
|P^+\gamma^+(f)|^2_{\atop L^2\left(\Gamma_+\right)}\leq {2\alpha\over
1-k_0^2}({\cal A}^\alpha f, f)_\Theta\leq C\alpha| f|^2_{\atop {\cal
A}^\alpha} 
\eeqna
where $k_0$ is such that
\quad$
\|{\cal K}P^+\|_{{\cal L}(L^2({\cal S}_+), L^2({\cal S}_+))}<k_0$.
 \end{lem}
\smallskip
\noindent
{\bf Proof.}
The outline of the proof, which is analogous of that of Lemma 4.4 in
\cite{DM} is as follows. We apply  Green's formula 
(\ref{LEMMAGREEN-1}) to $f=g$ times a  cutoff function $\rho_R
(|v|^2/2) =\rho(|v|^2/2R))$ with $\rho\in C^\infty (\R^+)$ such that
$0\leq \rho\leq 1,\,\,\, \rho(u)=1$ for $u<1$ and $\rho(u) =0$ for
$u>2$. We obtain, thanks to Lemma \ref{Nullspace} and hypothesis
\ref{HYPOPROJECT}:
\seqna
\label{LEMMATRACE2}
2\a({\cal A}^\a \rho_R, \rho_R f)_\Theta \geq
(1-k_0^2)
\int_{\Gamma_+}|v_x| P^+\gamma^+(f)|^2|\rho_R|^2d\Gamma.
\eeqna
Using the boundedness of $\phi$ leads to 
\seqna
\label{LEMMATRACE4}
2\a({\cal A}^\a \rho_R, \rho_R f)_\Theta \leq
|{\cal A}^\a f|_{L^2(\Theta)} + |f|_{L^2(\Theta)} +   {C\over R}
|f|_{L^2(\Theta)}.
\eeqna
Finally, taking the limit $R\to \infty$ in this estimate and recalling
(\ref{LEMMATRACE2}) allows us to conclude.

\smallskip
Next, we give some a priori estimates on the projections ${\cal
C}^\pm$. Note that, for $f\in D({\cal A}_\eta^\a)$, we have
$\gamma^-(f) ={\cal K}_\eta(\gamma^+(f))$. By Lemma
(\ref{ORTHOPROJECT1}) and using the form (\ref{PERTUBEDOPERATOR}) of ${\cal K}_{\eta}$, the
orthogonal projection of the trace on ${\cal C}^-$ reads
$Q^-\gamma^-(f) ={1\over 1+\eta} {\cal J}Q^+ \gamma^+(f)$.
Therefore,  there exists a function   $q(f) = q(x, \xi,|v|),\quad x=0, 1,\,\, \xi\in \R^2,\, |v|>0$,
such that
\seqna
\label{FUNCTION-Q}
q=(1+\eta)Q^-\gamma^-(f),\quad\mbox{on}\,\, \,\,\,\Gamma_-,\qquad q =
Q^+\gamma^+(f),  \quad\mbox{on}\,\, \,\,\,\Gamma_+.
\eeqna

\smallskip
We deduce 
\begin{lem} 
If $f\in D({\cal A}_\eta^\alpha)$, then there exists a constant $C>0$
such that
\seqna
|q(f)|^2_{\Gamma, R}\leq  C\left(\alpha |f|^2_{{\cal
A}^\alpha} +R|f|_{L^2(\Theta)}\right).
\eeqna
where, for $R>0$, $|\cdot|^2_{\Gamma, R}$ is defined by
$|\varphi|^2_{\Gamma, R} =\int_{\Gamma,\, |v|\leq R}|v_x||\varphi|^2
d\Gamma$.

\end{lem}
\smallskip
\noindent
{\bf Proof.}
The key of the proof is the control of   $Q^-\gamma^-(f) =
Q^+\gamma^+(f)$. Notice that, if $f\in D({\cal A}^\a)$, then we can
straightforwardly verify that $\mbox{sgn}(v_x)\zeta(x) f\in H_0({\cal
A}^\a)$. We multiply ${\cal A}^\a f$ by $\rho_{\atop R}\mbox{sgn}
(v_x)\zeta(x)f$ with $\zeta$ a function such that
$\zeta(1)=1,\,\, \zeta(0) =-1$ and integrate on  $\Gamma$. Thus we
apply Green's formula (\ref{LEMMAGREEN-1}) and take the limit
$R\to\infty$ just as in the proof of Lemma \ref{LEMMATRACE0}.
 The remainder of the proof follows now straightforward from
\cite{DM}, so we omit the detail.

\bigskip
The existence for the kinetic problem can now be stated. It is given
by  following  proposition, the proof  of which can be found in Ref.    \cite{DM}.
\begin{pro}
\label{PROPOEXISTENCE} 
Under hypothesis \ref{HYPOKERNEL} and  \ref{Hypotrans},  
 there exists a solution $f^\alpha$ to   problem
$(P)$, such that
$f^\alpha \in L^\infty(0, T; L^2(\Theta))$, ${\cal A} f^\alpha
\in  L^\infty(0, T; L^2(\Theta))$,   $P^+\gamma^+(f^\alpha) \in 
L^\infty(0, T; L^2(\Gamma_+))$, $Q^+\gamma^+(f^\alpha) \in 
L^\infty(0, T; L_R^2(\Gamma_+))$, for all $R>0$, and the boundary
condition is satisfied in the sense that
$\,\,\,\,\, P^-\gamma^-(f^\alpha)={\cal B} P^+
\gamma^+(f^\alpha);\quad  Q^-\gamma^-(f^\alpha)= 
Q^+\gamma^+(f^\alpha)$. Moreover, we have
\seqna
 \label{CCRR0}
\int_0^T |P^+\gamma^+(f^\alpha)(t)|_{L^2(\Gamma_+)}^2 dt &\leq&
C\a^2 |F_I|_{L^2(\Theta)}^2  .
\eeqna
 \end{pro}
\sect{The rigorous derivation of the SHE model}\label{RIGOROUS}

\quad 
In this section, we prove the convergence of the solutions of
the Boltzmann-Poisson equation towards solutions of the SHE model
(with coupled energies). The proof naturally falls   into five steps. First, we recall the $L^2$ estimate obtaining in the previous section; secondly, we prove the weak convergence of $f^\a$ and establish the continuity equation. Thirdly, we prove that the current actually converges. Fourthly, we pass to the limit in the nonlinear term and we finish the section by investigating some properties of the diffusivity of $\D$.

\subsection{L$^2$-estimates}

We summarize the L$^2$ estimates deduced from the previous  section. In the following estimates, $C$ is a bound for different constants depending only on the initial data.
 
\begin{lem}
\label{WEAKCONVERGE1}
The following estimates are satisfied by solution
$f^\a$ of the kinetic problem $(P)$
constructed in the previous section under hypotheses \ref{HYPOKERNEL}
and 
\ref{Hypotrans};
\seqna
\label{WEAK1}
|f^\a|_{C^0(0,T; L^2(\Theta))} &\leq& |F_I|_{
L^2(\Theta)}, \\
\label{WEAK11}
|{\cal A}^\a f^\a|_{C^0(0,T; L^2(\Theta))} &\leq& |{\cal A}^\a F_I|_{
L^2(\Theta)}, \\
\label{WEAK2}
\int_0^T |P^+\gamma^+(f^\a)|^2_{L^2(\Gamma_+)} \,dt &\leq& C
\a^2|F_I|_{ L^2(\Theta)}^2, \\
\label{WEAK3}
\int_0^T |P^-\gamma^-(f^\a)|^2_{L^2(\Gamma_-)} \,dt &\leq& C
\a^2|F_I|_{ L^2(\Theta)}^2, \\
\label{WEAK5}
\int_0^T |q(f^\a)|^2_{\Gamma, R} \,dt&\leq& C_R  |F_I|_{
H^\a({\cal A}^\a)}^2,
\eeqna
where $C$ denotes generic constants independent of $\a$ and of the data and  $q(f^\a)$ is defined  by (\ref{FUNCTION-Q}).

\end{lem}

As  a straight consequence of (\ref{WEAK1}), $(f^\a)_\a$ admits a
subsequence (still denoted by  $(f^\a))$ that converges to a function
$f^0$ in ${L^\infty(0,T; L^2(\R^3\times (0,\infty)})$ weak star as
$\a\to0$. From the inequalities (\ref{WEAK2}) and (\ref{WEAK3}) we deduce that the convergence of the traces $P^+\gamma^+(f^\a)$ (resp.  $P^-\gamma^-(f^\a)$) 
  to 0 in  ${L^2(0,T; L^2(\Gamma_+)})$ (resp. ${L^2(0,T; L^2(\Gamma_-)})$) strongly towards  zero.  Furthermore, using the diagonal extraction process, 
(\ref{WEAK5}) shows that  $q(f^\a)$ converges to a function $q(x,
\underline \xi,|v|, t)$ with $x=0, 1$ in $L^2(0, T; L^2(\gamma\times
B_R))$ weak star for any $R$, where $B_R$ is the ball centered at 0
and radius $R$ in velocity space by  properties of $H(div)$ spaces.
In order to give a precise meaning to the limits of traces, we notice
that the boundedness of ${\cal A}^\a f^\a$ in $H^{-1}(0, T;
L^2(\Theta))$  by (\ref{MainModel1}); this implies (according to 
\cite{B}, \cite{GR}) the boundedness of the sequence 
$(v_xf^\alpha)_\a|_\Gamma$ in $H^{-1}0, T; H^{1/2}(\gamma\times
B_R)$. Thank to this estimate, we obtain the  convergence in the
distributional sense of the traces of $f^\a$ on $\Gamma\times \R^3$
to the traces of $f^0$. Finally, the traces $\gamma^\pm(f^0)$ on
$\Gamma\times \R_\pm^3$ satisfy $P^-\gamma^-(f^0) =
P^+\gamma^+(f^0)=0,\quad Q^-\gamma^-(f^0) =Q^+\gamma^+(f^0) = q$,
so  that
\seqna
\label{WEAK6}
f^0|_{\Gamma}=q,
\eeqna 
is independent of the angular part of the velocity variable.

\subsection{Weak convergence of $f^\a$ and the continuity equation}

\quad From Green's formula (\ref{LEMMAGREEN-1}) we deduce  
the following  weak formulation of Eq. (\ref{MainModel1})-(\ref{Bord3}):

\begin{lem}
\label{WeakSolution1}
Let $f^\a$ be a solution to the kinetic problem
$(P)$  given by Proposition 
\ref{PROPOEXISTENCE}. Then, $\forall\, \psi\in C^1(\Theta'\times [0,
T], \R^2)$ (i.e. twice continuosly differentiable and with compact
support in $\R_{\underline \xi}^2 \times \R^{+*} \times [0, T[)$,
such that $\psi(T,\cdot, \cdot)=0$, we have
\seqna
\nonumber 
&&\int_0^T\!\!\int_\Theta f^\alpha\left(\alpha
{\d  \psi \over \d t} + (\underline v \cdot\nabla_{\underline \xi}
- {\underline E}^\alpha \nabla_{\underline v})\psi    +
{1\over\alpha}  \left( v_x{\d \over \d x}   - (\underline
v\times B)\cdot  \nabla_{\underline v}\right)\psi\right)dtd\theta\quad
\\ &&
\label{WeakSolution2}
+ \alpha\int_\Theta f_I\psi(0, \underline \xi, \underline
v)\, d\theta  \,\,=\,\,{1\over\alpha} \left(\int_0^T\!\!\int_\Gamma^+ 
|v_x| \gamma^+(f^\alpha) \gamma^+(\psi) -{\cal B}^*\gamma^-(\psi)
\,dtd\Gamma\right).
\eeqna
\end{lem}
\noindent
{\bf Proof.}
Equation (\ref{WeakSolution2}) is established by writing the weak
formulation of the perturbed
problem (\ref{PERTUBEDOPERATOR})-(\ref{EXISTENCESOLUTBord}) and
letting $\eta\to 0$ in this weak formulation.

\bigskip
As a consequence, we can prove that $f^0$ is a function of the total
energy only.

\begin{lem}
\label{ENERGY-II}
The limit function $f^0$ is a function of $(\underline \xi, |v|, t)$
only, $f^0 = f^0(\underline \xi, \frac12|v|^2, t)$.

\end{lem}
\noindent
{\bf Proof.}
Using (\ref{WeakSolution2}) with a test function $\psi$ such that
$\gamma^\pm(\psi)=0$, we get  
\seqna
\label{ENERGY-II1}
&&\a^2\int_0^T\!\!\int_\Theta f^\alpha  
{\d  \psi \over \d t} \,dtd\theta + \a \int_0^T\int_\Theta  f^\a
(\underline v \cdot\nabla_{\underline \xi}\psi
- {\underline E}^\alpha \nabla_{\underline v}\psi) \,dtd\theta \\
\nonumber
&& \int_0^T\!\!\int_\Theta f^\alpha \left( v_x{\d \over \d x} \psi -
(\underline v\times B)\cdot  \nabla_{\underline
v}\psi\right)dtd\theta + \int_\Theta f_I\psi|_{t=0} \,\, d\theta =0.
\eeqna
Hence, when $\a\to0$ in (\ref{ENERGY-II1}), using the fact that
$f^\a$ is bounded in $L^\infty(0, T; L^2(\Theta))$, we get 
\seqna
\label{ENERGY-II2}
\int_0^T\!\!\int_\Theta f^0 \left( v_x{\d \over \d x} \psi -
(\underline v\times B)\cdot  \nabla_{\underline
v}\psi\right)dtd\theta =0.
\eeqna
Eq. (\ref{ENERGY-II2}) supplemented with (\ref{WEAK6}), proves
that $f^0$ is a distributional solution of equation 
\seqna
\label{ENERGY-II3}
\left(v_x{\d \over \d x}   -
(\underline v\times B)\cdot  \nabla_{\underline
v} \right)\, f^0 =0,\qquad f^0|_{\Gamma}=q.
\eeqna
The remainder of the proof is a straightforward adaptation of Lemma
\ref{ENERGYFUNCTION} where we replace (\ref{FunctEner4E2}) by
$\gamma^-(f) =q$.

\bigskip
\quad
It remains to prove  that Eqs. (\ref{CONT})-(\ref{NEWINIT})
are satisfied in a weak sense. This is the object of the following
\begin{lem}
\label{CONTINUITYEQUATION}
 For any test function $\psi(\underline \xi, \varepsilon, t)\in
C^2(\Theta'\times [0, T])$ with support compactly supported in
$\R_{\underline \xi}^2\times\R^+\times [0, T[$, we have
\seqna
\label{CONTINUITYEQUATION1}
\nonumber
\int_0^T\!\!\!\int_{\R^2\times \R^{+*}} \left(
4\pi \sqrt{2\varepsilon}\, F^\a  {\d \psi\over \d t}+ \underline
J \cdot \left(\nabla_{\underline\xi}-{\underline E}{\d\over
\d\eps}\right)\psi\right) \,\,dtd\theta' \\
 + \int_{\R^2\times \R^{+*}}  4\pi \sqrt{2\varepsilon}\,
F_I\psi|_{t=0} \,\,d\theta' =0.
\eeqna
\end{lem}
\smallskip\noindent
{\bf Proof.} We apply  (\ref{LEMMAGREEN-1}) with $\psi$ as a test
function. Since $\psi$ is even with respect to $v$, it is easy to
check that the right-hand side in (\ref{WeakSolution2}) vanishes.
On the one hand, 
\seqna
\int_0^T\!\!\int_\Theta f^\alpha  {\d\psi \over \d t} \,\,dtd\theta
\to \int_0^T\!\!\int_\Theta F {\d  \psi \over \d t} \,\,dtd\theta
\eeqna
for $f^\a\to F$ in  $L^\infty(0, T; L^2(\Theta))$ weak-star and
$\d_t\psi\in  L^1(0, T; L^2(\Theta))$.  Moreover, on the grounds that  $\int_\Theta f_I\psi(0, \underline \xi, \underline
v)\, d\theta = \int_{\R^2\times \R^{+*}}  4\pi \sqrt{2\varepsilon} F_I\psi|_{t=0}
\,\,d\theta'$ and  $\left(v_x{\d \over \d x}   - (\underline
v\times B)\cdot  \nabla_{\underline v}\right)\psi=0$, the result
obviously follows from the coarea formula.
 The continuity equation  (\ref{CONT}) follows from
(\ref{CONTINUITYEQUATION1}) by taking the limit $\a\to0$.

\subsection{Equation for the current}

Our goal in this section is to prove that the approximate current
$\underline J^\a$ converges weakly towards a current $\underline J$ 
that satisfies (in a weak sense) the current equation (\ref{CURR})
(with
$F$ the limit of $(f^\a)_\a)$.  We are first concerned with the
consideration of the  most general auxiliary problem of which
(\ref{EQUATIONHOMOLOGIQUE}) is a particular case. Given a function
$g(x, \o)$, find $\chi(x, \o)$ such that:
\seqna
\label{Homo1}
\left(-v_x{\d  \over \d x} + (\underline v\times
B)\nabla_{\underline v}\right)  \chi(x,\omega)
&=&g,\qquad (x,\omega)\in [0,1]\times \S^2\\
\label{bound3}
\gamma^+( \chi)  &=&{\cal K}^* (\gamma^-(\chi))
\qquad (x,\omega)\in {\cal S}_+.
\eeqna
Once $\chi$ is determined, we deduce as a corollary, that there exist
functions $\chi_i(x,
\o;\underline \xi, \eps),\,\,\, (i=y, z)$, solutions of problem 
 (\ref{Homo1})-(\ref{bound3}) with the right-hand side $g=\o_i$,
unique up to additive functions of $\underline \xi$ and $|v|$. 

The existence of $\chi$ is provided by the following lemma which is
shown by adapting the proof of Lemma 5.4  in \cite{DM} in a fairly
straightforward way. 
 
\begin{lem}
\label{LEMMAHOMOLOGI}
For all $g\in L^2([0,1]\times \S^2)$, the
problem (\ref{Homo1})-(\ref{bound3}) has a solution if and only if
$g$ satisfies
\seqna
\int_0^1\!\!\!\int_{\S^2} g(x,\omega)dx d\omega =0. 
\eeqna
Furthermore, if this condition is satisfied, the
solution $\chi$ is unique under the condition
\seqna
\label{CONSTRAINT}
\int_0^1\!\!\!\int_{\S^2} \chi(x,\omega)dx d\omega =0
\eeqna
and  the set of solutions is the one-dimensional linear manifold
$\{\chi + F(\varepsilon),\,\, \mbox{with $F(\varepsilon)$
arbitrary}\}$.
\end{lem}

\begin{cor}
The function $g=\o_i,\quad (i=y, z)$ satisfies  the assumptions of Lemma
\ref{LEMMAHOMOLOGI}. 
\end{cor}
Therefore, the auxiliary function $\underline \chi$ defined by
(\ref{EQUATIONHOMOLOGIQUE}) exists and is unique under the
constraint (\ref{CONSTRAINT}).

 \bigskip
We require the following regularity assumptions on $\chi_i$:
 
\begin{hypo}
\label{HYPOREGULARITY}
(i) $\chi_i,\,\,\, (i=y, z)$, belongs to $L^2([0,1]\times \S^2)$ for
almost every $(\underline \xi, \eps)\in \R_{\underline \xi}^2 \times
\R_\eps^+$ and are $C^1$ bounded functions on $\Theta$ away from the
set $\{v_x=0\}$.

 (ii) The functions $\omega_i \chi_j(x, \omega; \xi, \varepsilon)$
belongs to $L^1([0,1]\times \S^2)$ and
$\displaystyle \int_0^1\!\!\!\int_{\S^2}\omega_i
\chi_jd xd\omega$ is a $C^1$ function of $(\underline
\xi, \varepsilon)\in \R_{\underline \xi}^2\times \R_\varepsilon^+$,
uniformly bounded on $\R_\xi^2\times [0, \infty[_\varepsilon$ and
tending to 0 as $\varepsilon\to0$.
\end{hypo}

\smallskip
In order to appreciate the importance of this hypothesis and, in
particular, to realize that this hypothesis is not empty, we refer
the reader to the example of isotropic scattering, where Ref. \cite{DM}
have computed explicitly the diffusion coefficient.

From  Hypothesis \ref{HYPOREGULARITY} (ii), we deduce that the   
diffusivity tensor (\ref{COEFFDIFFUSION1})  is defined and is  
$C^1$ function of $(\underline \xi, \varepsilon)\in \R_{\underline
\xi}^2\times \R_\varepsilon^+$.

\smallskip
We  are now in a position to establish the current equation.
Actually we shall prove that $J^\a$ has a finite limit.

\begin{lem} 
\label{CURRENTLIMIT0}
$J^\alpha$ is bounded in $L^2([0, T]\times \Theta)$. As $\a\to0$, $J^\a
\rightharpoonup J$ in the distributional sense and its limit
satisfies a weak form of the current equation.
More precisely:
For any test function  $\underline\psi =(\psi_y, \psi_z)$ in
$C^1(\Theta'\times [0, T],
\R^2)$ (i.e. twice continuosly differentiable and with compact
support in
$\R_{\underline \xi}^2 \times \R^{+*} \times [0, T[$),  we have
 \seqna
\label{CURRENTLIMIT3}
\int_0^T\!\!\!\int_{\R^2\times
\R^{+*}}\underline J^\alpha\cdot \underline \psi \, dtd\theta' 
\rightharpoonup 
\int_0^T\!\!\!\int_{\R^2\times \R^{+*}} F\left(
\nabla_{\underline\xi} -{\underline E}{\d \over \d
\varepsilon} \right)\cdot(\D^T\underline \psi)\, dtd\theta',  
\eeqna
as $\a\to0$,  where $\D^T$ denotes the transpose of $\D$.
The right-hand side of equation (\ref{CURRENTLIMIT3}) is the weak
form  of the current equation (\ref{CONT-CURR}).
\end{lem}
\noindent
{\bf Proof.} (i) Since $J^\a$ is a combination of $\gamma^+(J^\a)$ and
$\gamma^-(J^\a)$, to prove the boundedness of the current  in
$L^2([0, T]\times \Theta)$, it is enough to prove
that $\gamma^+(J^\a)$  and
$\gamma^-(J^\a)$ separately are bounded in this space. We prove
it for $\gamma^+(J^\a)$, the proof being similar for
$\gamma^-(J^\a)$. First, we note that
\seqna
\la v \gamma^+(f^\a\ra = \la v P^+ \gamma^+(f^\a)\ra,
\eeqna 
by symmetry. Then, using coarea formula and estimate
(\ref{WEAK2}) we get 
\seqna
\label{BOUNDCURRENT1}
|\gamma^+(J^\a)|_{L^2([0, T]\times \Theta)}^2 &=&
\int_0^T\!\!\!\!\int_0^\infty {1\over \a^2}  
\left|\la v P^+\gamma^+(f^\alpha)\ra\right|^2 d\eps dt\\
\nonumber
&\leq& 
{C\over \a^2} \int_0^T\!\!\!\!\int_0^\infty \left|\int_{\Gamma_+}
\omega P^+\gamma^+(f^\alpha) d\Gamma\right|^2 d\eps dt\\
\nonumber
&\leq&
{C\over \a^2} \int_0^T|P_{+}\gamma^+(f^\alpha)|_{L^2(\Gamma_+)}^2\,
\, dt   \leq  C |F_I|_{L^2([0, T]\times \Theta)}^2. 
\eeqna
The result follows straightforwardly.

(ii) Now, we choose  $\phi$ as  a test function $\phi (\underline
\xi, \varepsilon, t)=\sqrt{2\varepsilon} \underline \psi(\underline
\xi,  \varepsilon, t) \cdot \underline \chi (\underline \xi,
\omega; \varepsilon, t){\bf 1}_\rho(v_x)$ where ${\bf 1}_\rho(v_x)$ is
given by
\seqna
\label{TESTFUNCTION}
{\bf 1}_\rho(v_x)= \left\{
\begin{array}{cc}
\displaystyle  0\,\,\, \quad |v_x|\leq \rho\\
\displaystyle 1 \quad |v_x| \geq 2\rho.
\end{array} \right.
\eeqna
The hypothesis \ref{HYPOREGULARITY} makes it possible  to pass to
the limit $\rho\to0$ in  (\ref{TESTFUNCTION}).
Because of (\ref{Homo1})-(\ref{bound3}) we have for $i\in\{y,z\}$:
\seqna
\label{CURRENTLIMT1}
\left (v_x{\d \over \d x} -(\underline v \times B)\cdot
\nabla_{\underline v}\right)(\sqrt{2\varepsilon} \psi_i \chi_i) &=&
\sqrt{2\varepsilon} \psi_i   \left(v_x{\d \over \d x}
-(\underline v
\times B)\cdot
\nabla_{\underline v}\right) (\omega_i)\\
\nonumber
&=&-\sqrt{2\varepsilon}\psi_i \chi_i,\qquad \mbox{in}\,\,\, \Theta
\eeqna
and
\seqna
\label{CURRENTLIMT2}
\gamma^+(\sqrt{\varepsilon}\psi_i \chi_i)-{\cal K}^*
\gamma^-(\sqrt{\varepsilon}\psi_i \chi_i)
=\sqrt{\varepsilon}\psi_i[\gamma^+(\chi_i) -{\cal K}^*
\gamma^-(\chi_i)]=0,\qquad \mbox{on}\,\,\, \Gamma.
\eeqna
So that, using (\ref{CURRENTLIMT1})-(\ref{CURRENTLIMT2}) together
with coarera formula (\ref{COAREA2}),  we get
\seqna
\label{CURRENTLIMT3}
&&{1\over \alpha} \int_0^T\!\!\!\!\int_{\Theta} f^\alpha
\left (v_x{\d \over \d x} -(\underline v \times B)\cdot
\nabla_{\underline v}\right)(\sqrt{2\varepsilon} \underline\psi\,
\underline\chi)\,\, dtd\theta\\
\nonumber
&-&
{1\over \alpha} \int_0^T\!\!\!\!\int_{\Gamma^+}|v_x|f^\alpha(
\gamma^+(\sqrt{2\varepsilon}\cdot \underline \chi\, \underline\psi)-
{\cal K}^* \gamma^-(\sqrt{2\varepsilon})\cdot \underline
\chi\, \underline\psi)\,\,dtd\Gamma \\
\nonumber &=& -{1\over \alpha} \int_0^T\!\!\!\!\int_{\Theta}f^\alpha
\underline \omega\cdot \underline \psi\sqrt{2\varepsilon}dtd\theta
= - \int_0^T\int_{\Theta} \underline J^\alpha \cdot \underline
\psi dtd\theta' ,
\eeqna 
which justifies  the introduction of the
auxiliary function $\underline \chi$. Thus, the weak formulation
(\ref{WeakSolution2}) yields:
\seqna
\label{CURRENTLIMT4}
\nonumber
\int_0^T\!\!\!\!\int_{\Theta'}\underline J^\alpha \underline
\psi dtd\theta' &=&\alpha
\int_0^T\!\!\!\!\int_{\Theta}\sqrt{2\varepsilon}f^\alpha
\underline \chi \cdot{\d \over \d t} \underline \psi dtd\theta +
\alpha \int_{\Theta} \sqrt{2\varepsilon} f_I \underline \chi \cdot
\underline \psi|_{t=0} d\theta\\
&+&
\int_0^T\!\!\!\!\int_{\Theta}f^\alpha (\underline v \cdot
\nabla_{\underline
\xi} - \nabla_{\underline \xi} \phi \cdot \nabla_{\underline
v})(\sqrt{2\varepsilon}\underline \chi \cdot
\underline \psi)dtd\theta.
\eeqna
Now, by letting $\alpha\to0$ in (\ref{CURRENTLIMT4}), since
$\chi\in L^2([0,1]\times \S^2)$, for almost every $(\xi,
\varepsilon)\in
\R_{\underline \xi}^2\times
\R_\varepsilon^+$, the first and second terms on the right hand side
converge to 0. Then
\seqna
\label{CURRENTLIMT5}
\nonumber
\lim_{\alpha\to0}\int_0^T\!\!\!\!\int_{\Theta'}\underline
J^\alpha \underline \psi dtd\theta'&=&
 \int_0^T\!\!\!\!\int_{\Theta}f^\alpha
\underline v \cdot \nabla_{\underline \xi}
(\sqrt{2\varepsilon}\underline \chi \cdot\underline
\psi)dt d\theta  -   \int_0^T\!\!\!\!\int_{\Theta}
\nabla_{\underline \xi} \phi \cdot \nabla_{\underline
v})(\sqrt{2\varepsilon}\underline \chi \cdot \underline
\psi)dtd\theta \\
&=& L_1-L_2. 
\eeqna
On one hand, using (\ref{FunctEner2}) we have  
 \seqna
\nonumber
L_1 &=&
\int_0^T\!\!\!\int_{\Theta}f^0  \cdot
\nabla_{\underline \xi} ( {2\varepsilon}\underline \omega
(\underline \chi \cdot \psi))\,\,dt d\theta = 
\int_0^T\!\!\!\!\int_{\Theta}F\nabla_{\underline
 \xi}\cdot\left((2\varepsilon)^{3/2}\int_0^1\!\!\!\int_{\S^2}
\underline
\omega (\underline \chi   \underline \psi)dxd\omega\right) dt
d\theta \\
\label{CURRENTLIMT6}
 &=& \int_0^T\!\!\!\int_{\Theta'}F\nabla_{\underline \xi}
\cdot(\D^T \underline \psi)\,\, dtd\theta'.
\eeqna
On the  other hand,
using Lemma \ref{CURRENTDIFFERENTIAL2}, we get:
 \seqna
\label{CURRENTLIMT7} 
L_2 &=&
\int_0^T\!\!\! \int_{\Theta}f^0(\nabla_{\underline\xi}\phi 
\nabla_{\underline v})
(\sqrt{2\varepsilon}\underline \chi\cdot\underline \psi)\,\,
dtd\theta \\
\nonumber
&=&
\int_0^T\!\!\! \int_{\Theta'}F\nabla_{\underline\xi}\phi\cdot
\left(
\sqrt{2\varepsilon}\int_0^1\!\!\!\int_{\S^2} \nabla_{\underline
v} (\sqrt{2\varepsilon} \underline \chi \cdot \underline
\psi)dxd\omega\right) \,\, dt d\theta' \\
\nonumber
&=&
\int_0^T\!\!\!\int_{\Theta'}F \nabla_{\underline\xi}\phi\cdot {\d \over \d \varepsilon}
\left((2\varepsilon)^{3/2}\int_0^1\!\!\!\int_{\S^2}\underline
\omega(\underline \chi )dxd\omega \right) \,\,dt d\theta'= 
\int_0^T\!\!\!\int_{\Theta'}F \nabla_{\underline \xi} {\d \over \d
\varepsilon}(\D^T\underline \psi)\,\,dtd\theta';
\eeqna
so that, combining (\ref{CURRENTLIMT6}) and (\ref{CURRENTLIMT7}) to
conclude.

\subsection{Passing to the limit in the nonlinear term}

Notice that  the linear term of (\ref{WeakSolution2}) converges
because of the weak convergence of $f^\a$. To pass to the limit in
the nonlinear term of Eq. (\ref{WeakSolution2}), we need the strong
convergence of $E^\a$. For this task, we use the Poisson equation to
prove some compactness on the electric fields.

\smallskip
\quad
On the one hand, thanks to the classical regularizing properties of the
Laplacian, $F^\a (\underline \xi, \varepsilon, t)$ lies in
$L^\infty(0, T; L^2 (\Theta'))$, implies that,  since $\phi^\a$ solves
the Poisson equation,  $\phi^\a(\underline \xi , t)$ belongs to $L^\infty(0, T;
W^{2, p}(\R_{\underline\xi}^2))$, for all $\in[1, +\infty[$. In particular
$\nabla_{\underline \xi}\phi^\a\in L^\infty(0, T; W^{1, p}(\R_{\underline\xi}^2))$, for some $p>3$  and Sobolev's imbedding leads to $\nabla_{\underline \xi } \phi^\a\in L^\infty(\R_{\underline\xi}^2)$.

\quad 
On the other hand, by Poisson equation, we have 
\seqna
\label{POISSON-E1}
E^\a =\nabla_{\underline\xi}\Delta_{\underline\xi}^{-1}
(F^\a-C(\underline\xi))
\eeqna
The relation  (\ref{POISSON-E1}) together the continuity equation give  
\seqna
\label{POISSON-E2}
\d_tE^\a =\nabla_{\underline\xi}\Delta_{\underline\xi}^{-1}
\cdot\int f^\a \, v dxdv.
\eeqna
Hence, as ${\underline J}^\alpha$ is bounded in
$L^2(0,T;  L^2 (\Theta'))$ by  Lemma
\ref{CURRENTLIMIT0},  we conclude that $\d_t F^\alpha$ is bounded in
$L^2(0,T;  W^{-1, 2}(\Theta'))$ and
from  (\ref{POISSON-E2})  the regularizing properties of
the Laplacian now yield $\d_t \phi^\alpha$ bounded in
$L^\infty(0,T;  W^{1, 2}(\R_{\underline\xi}^2))$, which yields
$\d_t E^\alpha$  bounded in $L^\infty(0,T;
L^2(\R_{\underline\xi}^2))$.
 Thus, as  $W^{1, 2}(\R_{\underline\xi}^2)\hookrightarrow
L^2_{loc}(\R_{\underline\xi}^2)$, the Aubin-Lions Lemma    (see 
\cite{JPA}, \cite{L}) yields that the functional space
 \seqna \nonumber {\cal M} =\{E\in
L^\infty(0,T; W^{1, 2}(\R_{\underline\xi}^2)),\,\,\,\,  \d_t E  \in
L^\infty(0,T;  L^{2}(\R_{\underline\xi}^2))\},
\eeqna 
provided with the usual product norm, is compactly embedded
in $L^\infty(0,T;  L_{loc}^{2}(\R_{\underline\xi}^2))$ and
consequently, as $E^\alpha$ is bounded in ${\cal M}$,  there exists
a subsequence of $E^\alpha$ which converges strongly in
$L^\infty(0,T; L_{loc}^{2}(\R_{\underline\xi}^2))$.
Its remains to prove  that $E$  is the solution to the Poisson equation. To
do so, we notice that, $\phi^\a\in L^\infty(0, T; W^{2,
p}(\R_{\underline\xi}^2))$ and bounded in this space. Then, up an
extraction, we have

\smallskip
\centerline{$\phi^\a  \rightharpoonup  \phi\quad L^\infty(0, T; W^{2,
p}(\Theta'))\quad
\mbox{weak-star}.$}

\smallskip\noindent
Thus, it is easy to notice that $\phi$ is solution of Poisson equation
in the sense of distribution. 

\smallskip
Hence, we have

\begin{lem}
Under hypotheses \ref{HYPOCOUCHEINIT} and \ref{Hypotrans}, let
$(f^\a)$ be a family of a weak solution of the system
(\ref{MainModel1})-(\ref{MainModel2}). Then, extracting a
subsequence, the family $(E^\a)$ satisfies

\smallskip\smallskip
\centerline{$E^\a  \longrightarrow E\quad L^\infty(0, T; L^2_{loc}
(\R_{\underline\xi}^2)) \quad \mbox{strong}$,}

\noindent
as $\a\to0$.
\end{lem}

\subsection{Positivity of the diffusion tensor}\label{POSDIFFUSE}

We are now going to discuss some properties of the diffusivity of $\D$. We
show that $\D$ is a positive definite tensor. Therefore
(\ref{CONT})-(\ref{CURR}) are  well-posed. To this end, we prove

\begin{lem}
The diffusion tensor $\D$ defined by (\ref{COEFFDIFFUSION1}) is
positive definite; this means that, for all ${\underline\xi}\in \R^2$
and all $\eps_0>0$, there exists $C= C(\eps_0)>0$ such that:
\seqna
\label{Difften1}
(\D Y, Y) =\sum_{i,j}^2 \D_{ij}Y_iY_j\geq C|Y|^2 = C
\sum_{i=1}^2 Y_i^2,\qquad \forall\,  Y,{\underline\xi}\in
\R^2,\qquad \forall\,  \eps\geq \eps_0.
\eeqna 
\end{lem}
{\bf Proof.}  The proof is inspired from \cite{DM}. Let
$Y=(y_1, y_2)\in \R^2$ such that $|Y|>0$ and  
$\Phi(x,\omega)  = \sum_{i=1}^2  y_i \chi_i(x,\omega)$.
By virtue of (\ref{COEFFDIFFUSION1}),  we get
\seqna
\label{TENSOR1}
(\D Y,Y) &=& (2\eps)^{3/2}\int_0^1\int_{\S^2}
\left(\sum_{i=1}^2 \chi_i y_i\right)
\left(\sum_{i=1}^2 \omega_iy_i\right) dxd\omega
\eeqna
Inserting (\ref{EQUATIONHOMOLOGIQUE}) into  (\ref{TENSOR1}) yields
\seqna
\label{TENSOR2}
(\D Y, Y) &=& (2\eps)^{3/2}\int_0^1\int_{\S^2}\left(\sum_{i=1}^2
\chi_i y_i\right) \left(-v_x{\d \chi_i\over \d x} - (v\times
B)\nabla_{\underline v}) \chi_i y_i\right) dxd\omega\\
\nonumber
&=& (2\eps)^{3/2}({\cal A}^{0*}\Phi, \Phi )_{\cal S}.
\eeqna
On the other hand,  Green's formula leads to
\seqna
2 (\D Y, Y) &=&(2\eps)^{3/2} \left(\int_{{\cal S}_-} 
|\omega_x||\gamma^-(\Phi)|^2 d\omega  - \int_{ {\cal S}_+}
 |\omega_x||\gamma^+(\Phi)|^2 d\omega \right)\\
\nonumber
&\geq& (2\eps_0)^{3/2} \left(|\gamma^-(\Phi)|^2_{L^2({\cal
S}_-)} -|{\cal B^*} \gamma^-(\Phi)|^2_{L^2( {\cal S}_+)}\right)\geq0.
\eeqna
On the grounds that $|\omega_x|$ is uniformly bounded by positive constant, we
deduce that $\D$ is definite. Next,  assume that $\D Y\cdot
Y=0$. Then we deduce that $\Phi$ does not depend on $\underline
\omega$. So that ${\cal A}^{0*}\Phi={\underline \omega}\cdot Y$.
Since ${\cal A}^{0*}\Phi$ does not depend on $\underline \omega$, we
have $\underline \omega\cdot Y$ independent of $\underline \omega$ for
$\Omega$ in unit disk. This is only possible if $Y=0$
and  $\D$ is positive definite, and since  
$(\underline \xi,\eps)\longmapsto \D (\underline \xi,\eps)$ is
smooth, we deduce that $\D$ is uniformly positive definite. This
completes the proof.


\bigskip
\medskip\noindent
{\bf Acknowledgements.}
The author would like to express his
gratitude to P. Degond for bringing this problem
to his attention and T. Goudon  for their  fruitful comments and suggestions.   



\begin{thebibliography}{03}

\bibitem{Ar} A. A. Arsen'ev, Existence in the large of a weak
solution of Vlasov's system of equations, Z. Vycisl. Math. i Mat. Fiz
15 (1975) 136--147, 276 (Russian).

\bibitem{JPA} J.-P. Aubin, Un th\'eor\`eme de
compacit\'e, C. R. Acad. Sci. Paris 256 (1963), 5042-5044.

 
\bibitem{BBP}  H. Babovsky, C. Bardos, T. Platkwoski,
Diffusion approximation for a Knudsen gas in thin domain with
accommodation on the boundary, Asymptotic Analysis, 3, (1991), p.
265-289.

\bibitem{B}  C. Bardos,  Probl\`emes aux limites pour les
\'equations aux d\'eriv\'ees partielles du premier ordre \`a
coefficients r\'eels; Th\'eor\`emes d'approximations; application \`a
l'\'equation de transport, Ann. Sci. Ec. Norm. Sup,  4, (1970),
pp. 185-233. 

 
 
\bibitem{BSS} C. Bardos, C. Santos,   R. Sentis, Diffusion,
Approximation and Computation of
the Critical Size". Trans. Amer. Math. Soc., 284 A
(1984), 617--649.

\bibitem{BP}  R. Beals, V. Protopopescu,  Abstract Time Dependent Transport Problems,  J.  Math. Anal. Appl.   (1987), 121, 370.


\bibitem{BD} N. Benabdallah, P. Degond, On a hierarchy of
macrpscopic models for semiconductors, J. Maths. Phys. 37, (1996),
 3306--3333.

\bibitem{BLP}  A. Bensoussan,  J.L. Lions, G. C. Papanicolaou,
Boundary layers and homogenization of transport processes. J.
Publ. RIMS Kyoto Univ. 15, (1979),  53--157.

\bibitem{CIP}  C. Cercignani, R. Illner, M. Pulvirenti, The mathematical theory of dilute gases, Applied Mathematical Sciences, 106. Springer-Verlag, New  York, 1994. 


\bibitem{Ce} M. Cessenat, Th\'eor\`emes de traces $L^p$ pour les espaces de fonctions de la neutronique, C.R. Acad. sci. Paris. Ser. I, Math., 299 (1984), 
 831--834.


\bibitem{CH} S. Chapman,   T.G. Cowling, The mathematical
theory of non-uniform gases, Cambridge University Press, New-York,
1958.

\bibitem{DG}   J. Darroz\`es, J. P. Guiraud,  G\'en\'eralisation formelle du th\'eor\`eme H en pr\'esence de parois, Applications, C.R. Acad. sci. Paris. Ser. A, 262, (1966) 1368--1371.



\bibitem{DLMM}  P. Degond, V. Latocha, S. Mancini, A. Mellet, Diffusion dynamics of an electron gas between two plates, Methods Appl. Anal. 9 , (2002), 127--150. 
 
\bibitem{DM} P. Degond, S. Mancini, Diffusion driven by
collisions with the boundary. Asymptot. Anal. 27 (2001), no. 1,
47-73.

\bibitem{D2} C. Dogbe, Anomalous Diffusion Limit for a Knudsen
Gas in Thin Domain with   Accommodation on the Boundary.
J. of Stat.  Physics,  Vol. 100,   3/4, (2000), 603-632.

\bibitem{DSR} P. Dmitruk,  A. Sa\'ul,  L. Reyna,  High electric field approximation  to charge  transport  in semiconductor devices. Appl. Math. Lett. 5 (1992), no. 3, 99--102. 

 
\bibitem{F} H. Federer, Geometric  Measure Theory, Springer-Verlag,  New York,     1969.

\bibitem{GR} V. Girault, P-A. Raviart, Finite Element Methods for
Navier-Stokes Equations. Springer-Verlag, Berlin 1986.

\bibitem{GM}  T. Goudon, A. Mellet, On fluid limit for the semiconductors Boltzmann equation, J. Differential Equations,  J. Diff. Eq., 189,  (2003), 17-45.


 
\bibitem{Gy} Y. Guo,    Regularity for the Vlasov equations in a
half-space,   Indiana Univ. Math. J.  43  (1994),  no. 1, 255--320.
 
\bibitem{HJ} O. Hansen, A. J\"ungel, Analysis of Spherical harmonics
expansion model of physics, Math Models and Meth. in Appl. Sci. vol
14, n¡ 5, (2004) 759--774.

 

\bibitem{JP} A. J\"ungel, Y-J. Peng, A
hierarchy of hydrodynamic models for plasmas. Quasi-neutral limits
in the drift-diffusion equations, Asymptotic Analysis, 28 (2001), 49--73.

\bibitem{KR} M. G. Krein,  M. A. Rutman, Linear operator
leaving invariant a cone in Banach space. Uspehi Mat. Nauk. 3, no
1 (123) 3-95 (1948). Translation by  AMS Trans.  ser. 1, Vol 10,
199-325.

\bibitem{LK}   E. W. Larsen,   J. B Keller,  Asymptotic
Solutions of Neutron Transport Problem, J. Math. Phys., 15  (1974),
75--81. 



\bibitem{L} J. L. Lions, Quelques m\'ethodes de r\'esolution
de probl\`emes aux limites non lin\'eaires. Dunod,
Gauthier-Villars, 1969.

\bibitem{M}  S. Mischler, On the trace problem for solutions of Vlasov equation, Comm. in P.D.E., 25 (2000), 1415-1143.



\bibitem{Pa} K. Plaffelmoser, Global classical solution of the
Vlasov-Poisson system in three dimensions for general initial data, J.
Diff. Eq. 95 (1992), 281-303.

\bibitem{U} S. Ukai, Solutions of the Boltzmann equation, Patterns and Waves. Qualitative Analysis of Nonlinear Differential equations, North-Holland, Amsterdam, 1986,   37--96.



\bibitem{Z}  E. Zeidler, Nonlinear functional analysis and its applications. I. Fixed-point theorems. 
 Springer-Verlag, New York, 1986.

\end{thebibliography}
\end{document}